\documentclass [12pt]{article}

\usepackage {multicol}
\usepackage{amssymb,amsmath}
\usepackage{mathrsfs}
\usepackage{hyperref}
\usepackage{graphicx}
\usepackage{color}
\usepackage{caption}
\usepackage{esint}

\begin{document}

\title{Extension and Application of Deleting Items and Disturbing Mesh Theorem of Riemann Integral}

\author{\normalsize  Jingwei Liu\footnote{\scriptsize{Corresponding Author: Email: liujingwei03@tsinghua.org.cn (J.W Liu)}} \\
{\scriptsize{( School of Mathematics and System Sciences, Beihang University, Beijing,100083,P.R.China)}}
}

\date{\normalsize Sep. 09,2018}
\maketitle

\textbf{Abstract:}
The deleting items and disturbing mesh theorems of Riemann Integral are extended to multiple integral,line integral and surface integral respectively by constructing various of incomplete Riemann sum and non-Riemann sum sequences which converge to the same limit of classical Riemann sum. And, the deleting items and disturbing mesh formulae of Green's theorem, Stokes' theorem and divergence theorem (Gauss's or Ostrogradsky 's theorem) are also deduced. Then, the deleting items and disturbing mesh theorems of general Stokes' theorem on differential manifold are also derived.


\section{Introduction}
\label{}

 The deleting items theorem and disturbing mesh theorem of Riemann Integral  and the combination of deleting items and disturbing theorem of Riemann Integral are proposed for the one-dimensional Riemann Integral in [1]. We aim to extend the deleting item and disturbing mesh theorems to multiple integral,line integral and surface integral in this paper. Then, we develop the deleting items and disturbing mesh formulae for Green's theorem, Stokes' theorem and divergence theorem (Gauss's or Ostrogradsky 's theorem). At last, we derive the deleting items and disturbing mesh theorem of general Stokes' theorem in manifold analysis.

 As Green's theorem, Stokes' theorem and divergence theorem have strong physics and engineering background. Our researches not only discuss the integral based on Riemann sum in  mathematical analysis, but also show the uncertainty of physical measurement based on integral in the sense of limit theory especially in computation simulation with finite partition, and  reveal the complexity of mathematical modeling from the existences of numerous different physical formulas in finite status  approaching the same classical physical limit formulas and obeying same theory.

\section{Deleting Items and Disturbing Mesh Theorem of Multiple Integral}
\label{}

There is no uniform expression on multiple integral, the integral intervals can be $[a_1, b_1)\times [a_2, b_2)\times \cdots [a_n, b_n)$ , $[a_1, b_1]\times [a_2, b_2]\times \cdots [a_n, b_n]$, general intervals $C_1\times C_2\times \cdots C_n$ (where $C_k$ is bounded, unbounded, open, closed, half open in $\mathbb{R}^1$), or Jordan measurable set [2,3,4,5].  Definition 1 refers to [3].

\
\

\textbf{Definition 1} (Multiple Integral) Let $I=\{x=(x_1,x_2,\cdots,x_n)\in \mathbb{R}^n | a_i\leq x_i \leq b_i, i=1,\ldots,n \}$, $|I|:=\prod\limits_{i=1}^{n} (b_i-a_i)$. Let $f: I\mapsto \mathbb{R}$ be a real-valued and bounded function. Let intervals $I_1,\ldots,I_m$ be partition of $I$ induced by partitions of the coordinate intervals $[a_i,b_i],i=1,\ldots,n$, such that $I=\bigcup\limits_{k=1}^{m} I_k$, and no two of the intervals $I_1,\ldots,I_m$ have common interior points. Denote $P=\{I_1,\ldots,I_m \}$ , The mesh $\lambda(P):=\max\limits_{1\leq k \leq m} d(I_k)$ (the maximum among the diameters of the intervals of the partition P). For any distinguished points $\xi=\{\xi_1,\ldots,\xi_m \}$, $\xi_k \in I_k$, $k=1,\ldots,m$. The sum
\begin{equation}
\sigma(f,P,\xi):=\sum_{k=1}^{m} f(\xi_k) |I_k|=\sum_{k=1}^{m} f(\xi_k) m(I_k)
\end{equation}
is called  Riemann sum of $f$ corresponding to the partition of interval $I$ with distinguished points $(P,\xi)$, where $m(I_k)$ is the measure of $I_k$ .
\begin{equation}
\int_{I} f(x)dx:=\lim_{\lambda(P)\rightarrow 0} \sigma(f,P,\xi),
\end{equation}
Provides the limit exists, it is called the multiple Riemann integral of function $f$ over the interval $I$ and $f$ is integrable.
Or,
\begin{equation}
\displaystyle \int_{I} f(x)dx=\displaystyle \int\!\cdots\!\int_{I} f(x_1,x_2,\ldots,x_n) dx_1\cdots dx_n:=\lim_{\lambda(P)\rightarrow 0} \sum_{k=1}^{m} f(\xi_k) m(I_k).
\end{equation}
Particularly, if $m(I_1)=m(I_2)\cdots=m(I_m)$, $P$ is called equal partition, though $I_1, I_2\cdots I_m$ may not be congruent figures.

The deleting items and disturbing mesh theorems of Multiple Integral are as follows.

\hspace*{1cm}

\textbf{Theorem 1} (Deleting Items Theorem of Multiple Integral) Let $I=\{x=\{x_1,x_2,\cdots,x_n\}\in \mathbb{R}^n | a_i\leq x_i \leq b_i, i=1,\ldots,n \}$, $|I|:=\prod\limits_{i=1}^{n} (b_i-a_i)$. Let $f: I\mapsto \mathbb{R}$ be a real-valued and bounded function. Let intervals $I_1,\ldots,I_m$ be partition of $I$ induced by partitions of the coordinate intervals $[a_i,b_i],i=1,\ldots,n$, such that $I=\bigcup\limits_{k=1}^{m} I_k$, and no two of the intervals $I_1,\ldots,I_m$ have common interior points. Denote $P=\{I_1,\ldots,I_m \}$ , the mesh $\lambda(P):=\max\limits_{1\leq k \leq m} d(I_k)$ (the maximum among the diameters of the intervals of the partition P). For any distinguished points $\xi=\{\xi_1,\ldots,\xi_m \}$, $\xi_k \in I_k$, $k=1,\ldots,m$.  Denote $J = \{1,2, \cdots,m\}$. Then,

\begin{itemize}
\item[i)]  For a fixed natural number $K \in N^{+} (1 \leq K < m)$, denote $J_K = \{i_1 , \cdots ,i_K \} \subset J$, and $J\setminus J_{K}=\{1,2,\cdots,m\}\setminus \{i_1 , \cdots ,i_K \} $.\footnote{Given $m$, there are $C_{m}^{K}$ index sets of $J_K =\{i_1 , \cdots ,i_K \}$.}
Then,

\begin{equation}
\begin{array}{ll}
\displaystyle \lim\limits_{\lambda(P)\rightarrow 0} \sum\limits_{k\in J\setminus J_{K}} f(\xi_k) m(I_k)=\displaystyle \int_{I} f(x)dx.
\end{array}
\end{equation}
where $m(I_k)$ is the measure of $I_k$ .

\item[ii)] If $P$ is an equal partition, for any natural number $K(m)  (1 \leq K(m) < m)$, satisfying $\lim\limits_{m\rightarrow +\infty}\displaystyle \frac{K(m)}{m}=0 $, denote $J_{K(m)}=\{i_1 , \cdots ,i_{K(m)} \} \subset J$, and $J\setminus J_{K(m)}= \{1,2,\cdots,m\}\setminus \{i_1 , \cdots ,i_{K(m)} \} $. \footnote{Given $m$ and $K(m)$, there are $C_{m}^{K(m)}$ index sets of $J_{K(m)} =\{i_1 , \cdots ,i_{K(m)} \}$.}
Then,
\begin{equation}
\begin{array}{ll}
\displaystyle \lim\limits_{\lambda(P)\rightarrow 0} \sum\limits_{k\in J\setminus J_{{K(m)}}} f(\xi_k) m(I_k)=\displaystyle \int_{I} f(x)dx.
\end{array}
\end{equation}
where $m(I_k)$ is the measure of $I_k$ .

\end{itemize}

\textbf{Proof}  Since $f$ is bounded, $\exists M>0$, so that $|f|<M$.
As $f$ is integrable over $I$, $\forall \epsilon>0 $, $\exists \delta=\epsilon>0$, while $\lambda(P)<\delta$,

\begin{equation*}
|\sum\limits_{k=1}^{m} f(\xi_k) m(I_k)-\int_{I} f(x)dx|<\epsilon.
\end{equation*}

Since $\displaystyle \lim\limits_{\lambda(P)\rightarrow 0} \max\limits_{1\leq k \leq m} \{ m(I_k) \}=0$.
For above $\epsilon>0$, $\exists \delta=\epsilon>0$, while $\lambda(P)<\delta$,
\begin{equation*}
|\max\limits_{1\leq k \leq m} \{ m(I_k) \}|<\epsilon.
\end{equation*}

(1) Since
\begin{equation*}
\begin{array}{ll}
 &|\sum\limits_{k\in J\setminus J_{K}} f(\xi_k) m(I_k)-\int_{I} f(x)dx|\\
 =& |\sum\limits_{k=1}^{m} f(\xi_k) m(I_k)-\displaystyle \int_{I} f(x) dx - \sum\limits_{k\in  J_{K}} f(\xi_k) m(I_k)| \\
 \leq & |\sum\limits_{k=1}^{m} f(\xi_k) m(I_k)-\displaystyle \int_{I} f(x) dx|+ |\sum\limits_{k\in  J_{K}} f(\xi_k) m(I_k)| \\
 \leq & \varepsilon + KM\lambda(P)  \leq \varepsilon + KM\varepsilon  = (1+KM)\varepsilon \\
\end{array}
\end{equation*}
Then,
\begin{equation*}
\begin{array}{ll}
\displaystyle \lim\limits_{\lambda(P)\rightarrow 0} \sum\limits_{k\in J\setminus J_{K}} f(\xi_k) m(I_k)=\displaystyle \int_{I} f(x)dx
\end{array}
\end{equation*}

(2) As
\begin{equation*}
\begin{array}{ll}
 &|\sum\limits_{k\in J\setminus J_{{K(m)}}} f(\xi_k) m(I_k)-\int_{I} f(x)dx|\\
 =& |\sum\limits_{k=1}^{m} f(\xi_k) m(I_k)-\int_{I} f(x)dx + \sum\limits_{k\in  J_{{K(m)}}} f(\xi_k) m(I_k)| \\
 \leq & |\sum\limits_{k=1}^{m} f(\xi_k) m(I_k)-\int_{I} f(x)dx|+ |\sum\limits_{k\in  J_{{K(m)}}} f(\xi_k) m(I_k)| \\
 \leq & \varepsilon + M \sum\limits_{k\in  J_{{K(m)}}} m(I_k)\\
  =& \displaystyle \varepsilon + M m(I) \frac{ \sum\limits_{k\in  J_{{K(m)}}} m(I_k)}{m(I)} \\
\leq & \displaystyle  \varepsilon + M m(I) \frac{K(m)}{m} =(1+M m(I))\varepsilon \\
\end{array}
\end{equation*}
Then,
\begin{equation*}
\begin{array}{ll}
\displaystyle \lim\limits_{\lambda(P)\rightarrow 0} \sum\limits_{k\in J\setminus J_{{K(m)}}} f(\xi_k) m(I_k)=\displaystyle \int_{I} f(x)dx
\end{array}
\end{equation*}

\hfill $\Box$

\hspace*{1cm}

\textbf{Theorem 2} (Disturbing Mesh Theorem of Multiple Integral) Let $I=\{x=(x_1,x_2,\cdots,x_n)\in \mathbb{R}^n | a_i\leq x_i \leq b_i, i=1,\ldots,n \}$, $|I|:=\prod\limits_{i=1}^{n} (b_i-a_i)$. Let $f: I\mapsto \mathbb{R}$ be a real-valued and bounded function. Let intervals $I_1,\ldots,I_m$ be partition of $I$ induced by partitions of the coordinate intervals $[a_i,b_i],i=1,\ldots,n$, such that $I=\bigcup\limits_{k=1}^{m} I_k$, and no two of the intervals $I_1,\ldots,I_m$ have common interior points. Denote $P=\{I_1,\ldots,I_m \}$ , The mesh $\lambda(P):=\max\limits_{1\leq k \leq m} d(I_k)$ (the maximum among the diameters of the intervals of the partition P). For any distinguished points $\xi=\{\xi_1,\ldots,\xi_m \}$, $\xi_k \in I_k$, $k=1,\ldots,m$.

Let intervals $\tilde{I}_1,\ldots,\tilde{I}_m$ be distorted or inaccuracy grid of $I_1,\ldots,I_m$ respectively, and $\tilde{I}=\bigcup\limits_{k=1}^{m} \tilde{I}_k$, $\tilde{P}=\{\tilde{I}_1,\ldots,\tilde{I}_m \}$, $\lambda(\tilde{P}):=\max\limits_{1\leq k \leq m} d(\tilde{I}_k)$. If,
\begin{equation}
\begin{array}{ll}
\displaystyle \lim\limits_{m\rightarrow +\infty} \sum\limits_{k=1}^{m} m(I_k \triangle \tilde{I}_k)=0
\end{array}
\end{equation}
Then,
\begin{equation}
\begin{array}{ll}
\displaystyle \lim\limits_{\lambda(\tilde{P})\rightarrow 0} \sum\limits_{k=1}^{m}  f(\xi_k) m(\tilde{I}_k)=\displaystyle \int_{I} f(x)dx.
\end{array}
\end{equation}
where $m(I_k)$ is the measure of $I_k$ .

\textbf{Proof} Since
\begin{equation*}
\begin{array}{ll}
 \displaystyle|\sum\limits_{k=1}^{m}  f(\xi_k) m(\tilde{I}_k)-\int_{I} f(x) dx| \\
 =\displaystyle|\sum\limits_{k=1}^{m}  f(\xi_k) m(\tilde{I}_k)-\sum\limits_{k=1}^{m} f(\xi_k) m(I_k)+\sum\limits_{k=1}^{m} f(\xi_k) m(I_k)-\int_{I} f(x) dx|\\
 \leq \displaystyle|\sum\limits_{k=1}^{m}  f(\xi_k) m(\tilde{I}_k)-\sum\limits_{k=1}^{m} f(\xi_k) m(I_k)|+|\sum\limits_{k=1}^{m} f(\xi_k) m(I_k)-\int_{I} f(x) dx|\\
 = \displaystyle|\sum\limits_{k=1}^{m}  f(\xi_k) (m(\tilde{I}_k)- m(I_k))|+|\sum\limits_{k=1}^{m} f(\xi_k) m(I_k)-\int_{I} f(x) dx|\\
 \leq \displaystyle M \sum\limits_{k=1}^{m} (m(\tilde{I}_k \triangle I_k)) + |\sum\limits_{k=1}^{m} f(\xi_k) m(I_k)-\int_{I} f(x) dx|\\
\end{array}
\end{equation*}
and
\begin{equation*}
\begin{array}{ll}
\displaystyle \lim\limits_{m\rightarrow +\infty} \sum\limits_{k=1}^{m} m(I_k \triangle \tilde{I}_k)=0 , \  \ \lim\limits_{m\rightarrow +\infty} \sum\limits_{k=1}^{m} f(\xi_k) m(I_k)=\displaystyle \int_{I} f(x) dx
\end{array}
\end{equation*}

Then,
\begin{equation*}
\begin{array}{ll}
 \lim\limits_{m\rightarrow +\infty} \sum\limits_{k=1}^{m} f(\xi_k) m(\tilde{I}_k)=\displaystyle \int_{I} f(x) dx
\end{array}
\end{equation*}

This ends the proof.

\hfill $\Box$

Note that $f$ may not be defined on $\tilde{I}\backslash I$, though $f$ can be extended on it according to the values on intersect boundary. To avoid overelaborate discussion, we only consider  $\xi_k \in \tilde{I}\bigcap I$. For convenience, we assume that $f$ is defined on $\tilde{I}\bigcup I$ . And the distorted grid $\tilde{I}=\bigcup\limits_{k=1}^{m} \tilde{I}_k$ may not satisfy the two conditions: no two of the intervals have common interior points and $\bigcup\limits_{k=1}^{m} \tilde{I}_k=I$. In addition,
\begin{equation*}
\begin{array}{ll}
\displaystyle \lim\limits_{m\rightarrow +\infty} \sum\limits_{k=1}^{m} m(I_k \triangle \tilde{I}_k)=0
\end{array}
\end{equation*}
have concrete background. For example,\\

\textbf{Example 1.}
 $\displaystyle m(I_k)=\displaystyle \frac{|I|}{m}$, and $\displaystyle m(\tilde{I}_k)=\displaystyle \frac{(m^2-1)|I|}{m^2}$, $k=1,\cdots,m$. \\

\textbf{Example 2.}
$\displaystyle m(I_k)=\displaystyle \frac{|I|}{m}$, and $\displaystyle m(\tilde{I}_k)=\displaystyle \frac{(m^2+1)|I|}{m^2}$, $k=1,\cdots,m$.

\hspace*{1cm}

\textbf{Theorem 3} (Deleting Items and Disturbing Mesh Theorem of Multiple Integral) Let $I=\{x=\{x_1,x_2,\cdots,x_n\}\in \mathbb{R}^n | a_i\leq x_i \leq b_i, i=1,\ldots,n \}$, $|I|:=\prod\limits_{i=1}^{n} (b_i-a_i)$. Let $f: I\mapsto \mathbb{R}$ be a real-valued and bounded function. Let intervals $I_1,\ldots,I_m$ be partition of $I$ induced by partitions of the coordinate intervals $[a_i,b_i], i=1,\ldots,n$, such that $I=\bigcup\limits_{k=1}^{m} I_k$, and no two of the intervals $I_1,\ldots,I_m$ have common interior points. Denote $P=\{I_1,\ldots,I_m \}$ , The mesh $\lambda(P):=\max\limits_{1\leq k \leq m} d(I_k)$ (the maximum among the diameters of the intervals of the partition P). For any distinguished points $\xi=\{\xi_1,\ldots,\xi_m \}$, $\xi_k \in I_k$, $k=1,\ldots,m$.

Let intervals $\tilde{I}_1,\ldots,\tilde{I}_m$ be disturbed or inaccuracy grid of $I_1,\ldots,I_m$ respectively, and $\tilde{I}=\bigcup\limits_{k=1}^{m} \tilde{I}_k$, $\tilde{P}=\{\tilde{I}_1,\ldots,\tilde{I}_m \}$, $\lambda(\tilde{P}):=\max\limits_{1\leq k \leq m} d(\tilde{I}_k)$. If,
\begin{equation}
\begin{array}{ll}
\displaystyle \lim\limits_{m\rightarrow +\infty} \sum\limits_{k=1}^{m} m(I_k \triangle \tilde{I}_k)=0
\end{array}
\end{equation}
Denote $J = \{1,2, \cdots,m\}$.
Then,
\begin{itemize}
\item[i)]  For a fixed natural number $K \in N^{+} (1 \leq K < m)$, denote $J_K = \{i_1 , \cdots ,i_K \} \subset J$, and $J\setminus J_{K}=\{1,2,\cdots,m\}\setminus \{i_1 , \cdots ,i_K \} $.
Then,

\begin{equation}
\begin{array}{ll}
\displaystyle \lim\limits_{\lambda(\tilde{P})\rightarrow 0} \sum\limits_{k\in J\setminus J_{K}} f(\xi_k) m(\tilde{I}_k)=\displaystyle \int_{I} f(x) dx.
\end{array}
\end{equation}
where $m(I_k)$ is the measure of $I_k$ .

\item[ii)]If $P$ is an equal partition. For any natural number $K(m) (1 \leq K(m) < m)$, satisfying $\lim\limits_{m\rightarrow +\infty}\displaystyle \frac{K(m)}{m}=0 $, denote $J_{K(m)}=\{i_1 , \cdots ,i_{K(m)} \} \subset J$, and $J\setminus J_{K(m)}=\{1,2,\cdots,m\}\setminus \{i_1 , \cdots ,i_{K(m)} \} $.
Then,
\begin{equation}
\begin{array}{ll}
\displaystyle \lim\limits_{\lambda(\tilde{P})\rightarrow 0} \sum\limits_{k\in J\setminus J_{{K(m)}}} f(\xi_k) m(\tilde{I}_k)=\displaystyle \int_{I} f(x) dx.
\end{array}
\end{equation}
where $m(I_k)$ is the measure of $I_k$ .

\end{itemize}

\textbf{Proof}
(1) Since
\begin{equation*}
\begin{array}{ll}
\displaystyle \lim\limits_{m\rightarrow +\infty} \sum\limits_{k=1}^{m} m(I_k \triangle \tilde{I}_k)=0
\end{array}
\end{equation*}
Then,
\begin{equation*}
\begin{array}{ll}
\displaystyle \lim\limits_{\lambda(P)\rightarrow 0} \lambda(\tilde{P})=0
\end{array}
\end{equation*}
According to Theorem 1 and 2, we obtain
\begin{equation*}
\begin{array}{ll}
\lim\limits_{\lambda(\tilde{P})\rightarrow 0} \sum\limits_{k\in J\setminus J_{K}} f(\xi_k) m(\tilde{I}_k)
=\lim\limits_{\lambda(P)\rightarrow 0} \sum\limits_{k\in J\setminus J_{K}} f(\xi_k) m(I_k)\\
=\lim\limits_{\lambda(P)\rightarrow 0} \sum\limits_{k=1}^{m} f(\xi_k) m(I_k)
=\displaystyle \int_{I} f(x) dx
\end{array}
\end{equation*}

(2) According to Theorem 1 and 2, we obtain
\begin{equation*}
\begin{array}{ll}
\displaystyle \lim\limits_{\lambda(\tilde{P})\rightarrow 0} \sum\limits_{k\in J\setminus J_{{K(m)}}} f(\xi_k) m(\tilde{I}_k)
=\lim\limits_{\lambda(P)\rightarrow 0} \sum\limits_{k\in J\setminus J_{{K(m)}}} f(\xi_k) m(I_k)\\
=\lim\limits_{\lambda(P)\rightarrow 0} \sum\limits_{k=1}^{m} f(\xi_k) m(I_k)
=\displaystyle \int_{I} f(x) dx.
\end{array}
\end{equation*}

\hfill $\Box$

Note: Theorem 1,2,3 can be easily extended to the integral domain of Jordan-measurable set in $\mathbb{R}^n$.

\section{Deleting Items and Disturbing Mesh Theorem of Line Integral}
\label{}

There are two types of linear integral:  first type linear integral and second type linear integral, or scalar linear integral and vector linear integral [5,7].  Definition 2 and 3 refer to [7].

\textbf{Definition 2} (Scalar Line Integral) Let $\mathbf{x}:[a,b]\rightarrow \mathbb{R}^3 $ be a path of class $C^1$, and $f: D\subseteq \mathbb{R}^3 \rightarrow \mathbb{R}$ be a continuous function whose domain $D$ contains the image of $\mathbf{x}$ (so that the composite $f(\mathbf{x}(t))$ is defined). Let $a=t_0<t_1<\cdots<t_m=b$ be a partition of $[a,b]$. Let $t_k^{*}$ be an arbitrary points in the $k$-th subinterval $[t_{k-1},t_{k}]$ of the partition, $\triangle t_k = t_k - t_{k-1}$. Then we consider the Riemann sum
\begin{equation}
\sum_{k=1}^{m} f(\mathbf{x}(t_k^{*})) \triangle s_k,
\end{equation}
where  $\triangle s_k=\int_{t_{k-1}}^{t_k} \|\mathbf{x}^{'}(t) \|dt\approx \|\mathbf{x}^{'}(t_k^{*}) \|\triangle t_k $ is the length of the $k$-th segment of $\mathbf{x}$ (i.e. the portion of $\mathbf{x}$ defined for $t_{k-1}\leq t\leq t_{k}$ ). If the limit of Riemann sum exists as $\max_{\triangle t_k} \rightarrow 0$, we define the limit as the scalar line integral
\begin{equation}
\int_{\mathbf{x}} f ds=\displaystyle \lim_{\max\{\triangle t_k\} \rightarrow 0} \sum_{k=1}^{m}f(\mathbf{x}(t_k^{*})) \triangle s_k = \int_{a}^{b} f(\mathbf{x}(t)) \|\mathbf{x}'(t)\| dt.
\end{equation}
where $ds=\|\mathbf{x}'(t)\| dt$.

\textbf{Definition 3} (Vector Line Integral) Let $\mathbf{x}: [a,b]\rightarrow \mathbb{R}^n $ be a path of class $C^1$ ($n\geq 2$). Let $\mathbf{F}$ be a vector field defined on $D\subseteq \mathbb{R}^n$ such that $D$ contains the image of $\mathbf{x}$. Assume that $\mathbf{F}$ varies continuously along $\mathbf{x}$. Let $ a=t_0<t_1<\cdots<t_m=b$ be a partition of $[a,b]$. Let $t_k^{*}$ be an arbitrary points in the $k$-th subinterval $[t_{k-1},t_{k}]$ of the partition, $\triangle t_k = t_k - t_{k-1}$. Then we consider the Riemann sum
\begin{equation}
\sum_{k=1}^{m} \mathbf{F}(\mathbf{x}(t_k^{*})) \cdot \triangle \mathbf{s}_k,
\end{equation}
where  $\triangle \mathbf{s}_k= \mathbf{x}^{'}(t_k^{*})\triangle t_k $. If the limit of Riemann sum exists as $\max\{\triangle t_k\} \rightarrow 0$, we define the limit as vector line integral
\begin{equation}
\int_{\mathbf{x}} \mathbf{F} \cdot d \mathbf{s} =\displaystyle \lim_{\max\{\triangle t_k\} \rightarrow 0} \sum_{k=1}^{m} \mathbf{F}(\mathbf{x}(t_k^{*})) \cdot \triangle \mathbf{s}_k = \int_{a}^{b} \mathbf{F}(\mathbf{x}(t))\cdot \mathbf{x}'(t)dt.
\end{equation}
where $ d \mathbf{s} = \mathbf{x}^{'}(t)dt $.

The deleting items and disturbing mesh theorems of scalar line integral are as follows.

\hspace*{1cm}

\textbf{Theorem 4} (Deleting Items Theorem of Scalar Line Integral) Let $\mathbf{x}: [a,b]\rightarrow \mathbb{R}^3 $ be a path of class $C^1$. $f: D\subseteq \mathbb{R}^3 \rightarrow \mathbb{R}$ be a continuous function whose domain $D$ contains the image of $\mathbf{x}$ ,and the scalar line integral of $f$ over $\mathbf{x}$ exists.
Let $a=t_0<t_1<\cdots<t_m=b$ be a partition of $[a,b]$, $I_k=[t_{k-1},t_{k}]$, $\triangle t_k = t_k - t_{k-1}$. Let $t_k^{*}$ be an arbitrary points in $I_k$. $P=\{I_1,\cdots,I_m\} $.
Denote $J = \{1,2, \cdots,m\}$. Then,
\begin{itemize}
\item[i)]  For a fixed natural number $K \in N^{+} (1 \leq K < m)$, denote $J_K = \{i_1 , \cdots ,i_K \} \subset J$, and $J\setminus J_{K}=\{1,2,\cdots,m\}\setminus \{i_1 , \cdots ,i_K \} $.
Then,
\begin{equation}
\begin{array}{ll}
\displaystyle \lim_{\max\{\triangle t_k\} \rightarrow 0} \sum\limits_{k\in J\setminus J_{K}} f(\mathbf{x}(t_k^{*})) \triangle s_k= \int_{\mathbf{x}} f ds.
\end{array}
\end{equation}

\item[ii)] If $P$ is an equal partition, for any natural number $K(m) (1 \leq K(m) < m)$, satisfying $\lim\limits_{m\rightarrow +\infty}\displaystyle \frac{K(m)}{m}=0 $, denote $J_{K(m)}=\{i_1 , \cdots ,i_{K(m)} \} \subset J$, and $J\setminus J_{K(m)}=\{1,2,\cdots,m\}\setminus \{i_1 , \cdots ,i_{K(m)} \} $.
Then,
\begin{equation}
\begin{array}{ll}
\displaystyle \lim_{\max\{\triangle t_k\} \rightarrow 0} \sum\limits_{k\in J\setminus J_{K(m)}} f(\mathbf{x}(t_k^{*})) \triangle s_k= \int_{\mathbf{x}} f ds.
\end{array}
\end{equation}

\end{itemize}
where $\triangle s_k=\int_{t_{k-1}}^{t_k} \|\mathbf{x}^{'}(t) \|dt\approx \|\mathbf{x}^{'}(t_k^{*}) \|\triangle t_k $, $ds=\|\mathbf{x}'(t)\| dt$.

\hspace*{1cm}

\textbf{Theorem 5} (Disturbing Mesh Theorem of Scalar Line Integral) Let $\mathbf{x}: [a,b]\rightarrow \mathbb{R}^3 $ be a path of class $C^1$. $f: D \subseteq \mathbb{R}^3 \rightarrow \mathbb{R}$ be a continuous function whose domain $D$ contains the image of $\mathbf{x}$, and the scalar line integral of $f$ over $\mathbf{x}$ exists.

Let $a=t_0<t_1<\cdots<t_m=b$ be a partition of $[a,b]$, $I_k=[t_{k-1},t_{k}]$, $\triangle t_k = t_k - t_{k-1}$. Let $t_k^{*}$ be an arbitrary points in $I_k$. $P=\{I_1,\cdots,I_m\} $.

Let $a=\tilde{t}_0<\tilde{t}_1<\cdots<\tilde{t}_m=b$ be a distorted partition of $a=t_0<t_1<\cdots<t_m=b$, $\tilde{I}_k=[\tilde{t}_{k-1},\tilde{t}_{k}]$, $\triangle \tilde{t}_k = \tilde{t}_k - \tilde{t}_{k-1}$, $\tilde{P}=\{\tilde{I}_1,\cdots,\tilde{I}_m\}$, satisfying
\begin{equation*}
\displaystyle \lim\limits_{m\rightarrow +\infty} \sum\limits_{k=1}^{m} m(I_k \triangle \tilde{I}_k)=0
\end{equation*}
Then,
\begin{equation}
\begin{array}{ll}
\displaystyle \lim_{\max\{\triangle \tilde{t}_k\} \rightarrow 0} \sum\limits_{k=1}^m  f(\mathbf{x}(t_k^{*})) \triangle \tilde{s}_k= \int_{\mathbf{x}} f ds.
\end{array}
\end{equation}
where  $\triangle \tilde{s}_k\approx \|\mathbf{x}^{'}(t_k^{*}) \|\triangle \tilde{t}_k $, $ds=\|\mathbf{x}'(t)\| dt$.

\hspace*{1cm}

\textbf{Theorem 6} (Deleting Items and Disturbing Mesh Theorem of Scalar Line Integral) Let $\mathbf{x}: [a,b]\rightarrow \mathbb{R}^3 $ be a path of class $C^1$. $f: D\subseteq \mathbb{R}^3 \rightarrow \mathbb{R}$ be a continuous function whose domain $D$ contains the image of $\mathbf{x}$, and the scalar line integral of $f$ over $\mathbf{x}$ exists.

Let $a=t_0<t_1<\cdots<t_m=b$ be a partition of $[a,b]$, $I_k=[t_{k-1},t_{k}]$, $\triangle t_k = t_k - t_{k-1}$. Let $t_k^{*}$ be an arbitrary points in $I_k$. $P=\{I_1,\cdots,I_m\} $.

Let $a=\tilde{t}_0<\tilde{t}_1<\cdots<\tilde{t}_m=b$ be a distorted partition of $a=t_0<t_1<\cdots<t_m=b$, $\tilde{I}_k=[\tilde{t}_{k-1},\tilde{t}_{k}]$, $\triangle \tilde{t}_k = \tilde{t}_k - \tilde{t}_{k-1}$, $\tilde{P}=\{\tilde{I}_1,\cdots,\tilde{I}_m\}$, satisfying
\begin{equation*}
\displaystyle \lim\limits_{m\rightarrow +\infty} \sum\limits_{k=1}^{m} m(I_k \triangle \tilde{I}_k)=0
\end{equation*}
Denote $J = \{1,2, \cdots,m\}$.Then,
\begin{itemize}
\item[i)]  For a fixed natural number $K \in N^{+} (1 \leq K < m)$, denote $J_K = \{i_1 , \cdots ,i_K \} \subset J$, and $J\setminus J_{K}=\{1,2,\cdots,m\}\setminus \{i_1 , \cdots ,i_K \} $.
Then,
\begin{equation}
\begin{array}{ll}
\displaystyle \lim_{\max\{\triangle \tilde{t}_k\} \rightarrow 0} \sum\limits_{k\in J\setminus J_{K}} f(\mathbf{x}(t_k^{*})) \triangle \tilde{s}_k= \int_{\mathbf{x}} f ds.
\end{array}
\end{equation}

\item[ii)]If $P$ is an equal partition. For any natural number $K(m) (1 \leq K(m) < m)$, satisfying $\lim\limits_{m\rightarrow +\infty}\displaystyle \frac{K(m)}{m}=0 $, denote $J_{K(m)}=\{i_1 , \cdots ,i_{K(m)} \} \subset J$, and $J\setminus J_{K(m)}=\{1,2,\cdots,m\}\setminus \{i_1 , \cdots ,i_{K(m)} \} $.
Then,
\begin{equation}
\begin{array}{ll}
\displaystyle \lim_{\max\{\triangle \tilde{t}_k\} \rightarrow 0} \sum\limits_{k\in J\setminus J_{K(m)}} f(\mathbf{x}(t_k^{*})) \triangle \tilde{s}_k= \int_{\mathbf{x}} f ds.
\end{array}
\end{equation}
\end{itemize}
where  $\triangle s_k\approx \|\mathbf{x}^{'}(t_k^{*}) \|\triangle t_k $, $\triangle \tilde{s}_k\approx \|\mathbf{x}^{'}(t_k^{*}) \|\triangle \tilde{t}_k $, $ds=\|\mathbf{x}'(t)\| dt$.


The deleting items and disturbing mesh theorems of  vector line integral are as follows.

\hspace*{1cm}

\textbf{Theorem 7} (Deleting Items Theorem of Vector Line Integral) Let $\mathbf{x}:[a,b]\rightarrow \mathbb{R}^n $ be a path of class $C^1$ ($n\geq 2$). Let $\mathbf{F}$ be a vector field defined on $D\subseteq \mathbb{R}$ such that $D$ contains the image of $\mathbf{x}$. Assume that $\mathbf{F}$ varies continuously along $\mathbf{x}$.

Let $a=t_0<t_1<\cdots<t_m=b$ be a partition of $[a,b]$, $I_k=[t_{k-1},t_{k}]$, $\triangle t_k = t_k - t_{k-1}$. Let $t_k^{*}$ be an arbitrary points in $I_k$. $P=\{I_1,\cdots,I_m\} $.
Denote $J = \{1,2, \cdots,m\}$. Then,

\begin{itemize}
\item[i)]  For a fixed natural number $K \in N^{+} (1 \leq K < m)$, denote $J_K = \{i_1 , \cdots ,i_K \} \subset J$, and $J\setminus J_{K}=\{1,2,\cdots,m\}\setminus \{i_1 , \cdots ,i_K \} $.
Then,

\begin{equation}
\begin{array}{ll}
\displaystyle \lim_{\max\{\triangle t_k\} \rightarrow 0} \sum\limits_{k\in J\setminus J_{K}} \mathbf{F}(\mathbf{x}(t_k^{*})) \cdot \triangle \mathbf{s}_k =\displaystyle \int_{\mathbf{x}} \mathbf{F} \cdot d \mathbf{s}.
\end{array}
\end{equation}

\item[ii)] If $P$ is an equal partition, for any natural number $K(m) (1 \leq K(m) < m)$, satisfying $\lim\limits_{m\rightarrow +\infty}\displaystyle \frac{K(m)}{m}=0 $, denote $J_{K(m)}=\{i_1 , \cdots ,i_{K(m)} \} \subset J$, and $J\setminus J_{K(m)}=\{1,2,\cdots,m\}\setminus \{i_1 , \cdots ,i_{K(m)} \} $.
Then,
\begin{equation}
\begin{array}{ll}
\displaystyle \lim_{\max\{\triangle t_k\} \rightarrow 0} \sum\limits_{k\in J\setminus J_{K(m)}} \mathbf{F}(\mathbf{x}(t_k^{*})) \cdot \triangle \mathbf{s}_k =\displaystyle \int_{\mathbf{x}} \mathbf{F} \cdot d \mathbf{s}.
\end{array}
\end{equation}
\end{itemize}
where $ d \mathbf{s} = \mathbf{x}^{'}(t)dt $, $\triangle \mathbf{s}_k= \mathbf{x}^{'}(t_k^{*})\triangle t_k $.

\hspace*{1cm}

\textbf{Theorem 8} (Disturbing Mesh Theorem of Vector Line Integral) Let $\mathbf{x}:[a,b]\rightarrow \mathbb{R}^n $ be a path of class $C^1$ ($n\geq 2$). Let $\mathbf{F}$ be a vector field defined on $D\subseteq \mathbb{R}$ such that $D$ contains the image of $\mathbf{x}$. Assume that $\mathbf{F}$ varies continuously along $\mathbf{x}$,and the vector line integral of $\mathbf{F}$ over $\mathbf{x}$ exists.

Let $a=t_0<t_1<\cdots<t_m=b$ be a partition of $[a,b]$, $I_k=[t_{k-1},t_{k}]$, $\triangle t_k = t_k - t_{k-1}$. Let $t_k^{*}$ be an arbitrary points in $I_k$. $P=\{I_1,\cdots,I_m\} $.

Let $a=\tilde{t}_0<\tilde{t}_1<\cdots<\tilde{t}_m=b$ be a distorted partition of $a=t_0<t_1<\cdots<t_m=b$, $\tilde{I}_k=[\tilde{t}_{k-1},\tilde{t}_{k}]$, $\triangle \tilde{t}_k = \tilde{t}_k - \tilde{t}_{k-1}$,  $\tilde{P}=\{\tilde{I}_1,\cdots,\tilde{I}_m\}$, satisfying
\begin{equation*}
\displaystyle \lim\limits_{m\rightarrow +\infty} \sum\limits_{k=1}^{m} m(I_k \triangle \tilde{I}_k)=0
\end{equation*}
Then,
\begin{equation}
\begin{array}{ll}
\displaystyle \lim_{\max\{\triangle \tilde{t}_k\} \rightarrow 0} \sum\limits_{k=1}^m \mathbf{F}(\mathbf{x}(t_k^{*})) \cdot \triangle \tilde{\mathbf{s}}_k =\displaystyle \int_{\mathbf{x}} \mathbf{F} \cdot d \mathbf{s}.
\end{array}
\end{equation}
where $ d \mathbf{s} = \mathbf{x}^{'}(t)dt $, $\triangle \tilde{\mathbf{s}}_k= \mathbf{x}^{'}(t_k^{*})\triangle \tilde{t}_k $.

\hspace*{1cm}

\textbf{Theorem 9} (Deleting Items and Disturbing Mesh Theorem of Vector Line Integral) Let $\mathbf{x}: [a,b]\rightarrow \mathbb{R}^n $ be a path of class $C^1$ ($n\geq 2$). Let $\mathbf{F}$ be a vector field defined on $D\subseteq \mathbb{R}$ such that $D$ contains the image of $\mathbf{x}$. Assume that $\mathbf{F}$ varies continuously along $\mathbf{x}$, and the vector line integral of $\mathbf{F}$ over $\mathbf{x}$ exists.

Let $a=t_0<t_1<\cdots<t_m=b$ be a partition of $[a,b]$. Let $t_k^{*}$ be an arbitrary points in the $k$-th subinterval $I_k=[t_{k-1},t_{k}]$ of the partition, $\triangle t_k = t_k - t_{k-1}$.  $P=\{I_1,\cdots,I_m\} $.

Let  $a=\tilde{t}_0<\tilde{t}_1<\cdots<\tilde{t}_m=b$ be a distorted partition of $a=t_0<t_1<\cdots<t_m=b$, $\tilde{I}_k=[\tilde{t}_{k-1},\tilde{t}_{k}]$, $\triangle \tilde{t}_k = \tilde{t}_k - \tilde{t}_{k-1}$, $\tilde{P}=\{\tilde{I}_1,\cdots,\tilde{I}_m\}$, satisfying
\begin{equation*}
\displaystyle \lim\limits_{m\rightarrow +\infty} \sum\limits_{k=1}^{m} m(I_k \triangle \tilde{I}_k)=0
\end{equation*}
Denote $J = \{1,2, \cdots,m\}$.Then,
\begin{itemize}
\item[i)]  For a fixed natural number $K \in N^{+} (1 \leq K < m)$, denote $J_K = \{i_1 , \cdots ,i_K \} \subset J$, and $J\setminus J_{K}=\{1,2,\cdots,m\}\setminus \{i_1 , \cdots ,i_K \} $.
Then,

\begin{equation}
\begin{array}{ll}
\displaystyle \lim_{\max\{\triangle \tilde{t}_k\} \rightarrow 0} \sum\limits_{k\in J\setminus J_{K}} \mathbf{F}(\mathbf{x}(t_k^{*})) \cdot \triangle \tilde{\mathbf{s}}_k =\displaystyle \int_{\mathbf{x}} \mathbf{F} \cdot d \mathbf{s}.
\end{array}
\end{equation}

\item[ii)]If $P$ is an equal partition, for any natural number $K(m) (1 \leq K(m) < m)$, satisfying $\lim\limits_{m\rightarrow +\infty}\displaystyle \frac{K(m)}{m}=0 $, denote $J_{K(m)}=\{i_1 , \cdots ,i_{K(m)} \} \subset J$, and $J\setminus J_{K(m)}=\{1,2,\cdots,m\}\setminus \{i_1 , \cdots ,i_{K(m)} \} $.
Then,
\begin{equation}
\begin{array}{ll}
\displaystyle \lim_{\max\{\triangle \tilde{t}_k\} \rightarrow 0} \sum\limits_{k\in J\setminus J_{K(m)}} \mathbf{F}(\mathbf{x}(t_k^{*})) \cdot \triangle \tilde{\mathbf{s}}_k =\displaystyle \int_{\mathbf{x}} \mathbf{F} \cdot d \mathbf{s}.
\end{array}
\end{equation}
\end{itemize}
where $d \mathbf{s} = \mathbf{x}^{'}(t)dt $, $\triangle \mathbf{s}_k= \mathbf{x}^{'}(t_k^{*})\triangle t_k $ , $\triangle \tilde{\mathbf{s}}_k= \mathbf{x}^{'}(t_k^{*})\triangle \tilde{t}_k $.

The proofs of Theorem 4,5,6,7,8,9 are similar to Theorem 1,2,3 correspondingly, we omit the detail procedures.
And, Theorem 4,5,6,7,8,9 can also be discussed as [1] while variable of integration is $t$.

\section{Deleting Items and Disturbing Mesh Theorem of Surface Integrals}
\label{}

There are two types of surface integral: surface integral over scale field and surface integral over vector field [7]. They are also called  surface integral of the first type and surface integral of the second type [3,5]. Definition 4,5 refer to [7].

\textbf{Definition 4}. (Scaler Surface Integral) Let $S$ be a smooth parameterized surface, that is $S$ is the image of the $C^1$ map $\mathbf{X}: D\rightarrow \mathbb{R}^3$ , where $D$ is a connected ,bounded region in $\mathbb{R}^2$. Let $f$ be a continuous function defined on $S=\mathbf{X}(D)$.
Suppose $S$ is partitioned into finitely many small pieces $S_1,S_2,\cdots,S_m$, every two of them can only intersect along some parts of their boundaries, and the area of the $k$-th pieces is $\triangle S_k$. Let $c_k$ denote an arbitrary ``test point'' in the $k$-th piece. Then the scalar surface integral of $f$ over $S$ should be
\begin{equation}
\displaystyle \iint_{S} f dS=\displaystyle \lim_{\max\{\triangle S_k\} \rightarrow 0} \sum_{k=1}^{m} f(c_k) \triangle S_k.
\end{equation}
provided the limit exists.

Equally, for $S=\mathbf{X}(D)$, where $\mathbf{X}(u,v)=\{x(u,v),y(u,v),z(u,v)\}$, $(u,v)\in D$,  let $D_1,\cdots,D_m$ be any partition of $D$, the area of the $k$-th pieces be $\triangle D_k$, such that $S_k=\mathbf{X}(D_k)$, $k=1,\cdots,m$, $S_1, \cdots, S_m$ is a partition of $S$ and the area of the $S_k$ is $\triangle S_k$. Let $\xi_k$ denote an arbitrary ``test point'' in $D_k$, such that $c_k=\mathbf{X}(\xi_k)$,
$\displaystyle \triangle S_k \approx \|\frac{\partial X}{\partial u} \times \frac{\partial X}{\partial v} \|_{\xi_k} \triangle D_k $.
Then,
\begin{equation}
\begin{array}{ll}
\displaystyle \iint_{S} f dS=\displaystyle \lim_{\max\{\triangle S_k\} \rightarrow 0} \sum_{k=1}^{m} f(c_k) \triangle S_k  \\
=\displaystyle \lim_{\max\{\triangle D_k\} \rightarrow 0} \sum_{k=1}^{m} f(\mathbf{X}(\xi_k))\|\frac{\partial X}{\partial u} \times \frac{\partial X}{\partial v} \|_{\xi_k} \triangle D_k.\\
=\displaystyle \iint_{D} f(\mathbf{X}(u,v))\|\frac{\partial \mathbf{X}}{\partial u} \times \frac{\partial \mathbf{X}}{\partial v} \| dudv.
\end{array}
\end{equation}

\textbf{Definition 5}. (Vector Surface Integral) Let $\mathbf{X}: D\rightarrow \mathbb{R}^3$ be a smooth parameterized surface, where $D$ is a bounded region in the plane, and let $\mathbf{F}(x,y,z)$ be a continuous vector field whose domain includes $S=\mathbf{X}(D)$.

Suppose $S$ is partitioned into finitely many small pieces $S_1,S_2,\cdots,S_m$, every two of them can only intersect along some parts of their boundaries, and the area of the $S_k$ is $\triangle S_k$. Let $c_k$ denote an arbitrary ``test point'' in the $S_k$ piece.
Then the vector surface integral of $\mathbf{F}$ along $\mathbf{X}$ is denoted  as
\begin{equation}
\begin{array}{ll}
\displaystyle \iint_{\mathbf{X}} \mathbf{F}\cdot d\mathbf{S}
:=\displaystyle \lim_{\max\{\triangle S_k\} \rightarrow 0} \sum_{k=1}^{m} \mathbf{F}(c_k)\cdot \triangle \mathbf{S}_k
\end{array}
\end{equation}

Suppose $D_1,D_2,\cdots,D_m$ is a partition of $D$, such that $S_k=\mathbf{X}(D_k)$. The area of the $D_k$ is $\triangle D_k$. Let $\xi_k$ denote an arbitrary ``test point'' in $D_k$, such that $c_k=\mathbf{X}(\xi_k)$. Then
\begin{equation}
\begin{array}{ll}
\displaystyle \iint_{\mathbf{X}} \mathbf{F}\cdot d\mathbf{S}
=\displaystyle \lim_{\max\{\triangle D_k\} \rightarrow 0} \sum_{k=1}^{m} \mathbf{F}(\mathbf{X}(\xi_k))\cdot \mathbf{N}(\xi_k) \triangle \mathbf{D}_k\\
=\displaystyle \iint_{D} \mathbf{F}(\mathbf{X}(u,v))\cdot \mathbf{N}(u,v) du dv
\end{array}
\end{equation}
where $\mathbf{N}(u,v)=\displaystyle \frac{\partial \mathbf{X}}{\partial u} \times \frac{\partial \mathbf{X}}{\partial v}$.

The deleting items and disturbing mesh theorems of scaler surface integral are as follows.

\hspace*{1cm}

\textbf{Theorem 10} (Deleting Items Theorem of Scaler Surface Integral) Let $S$ be a smooth parameterized surface, that is $S$ is the image of the $C^1$ map $\mathbf{X}: D\rightarrow \mathbb{R}^3$,$(u,v)\in D$ , where $D$ is a connected, bounded region in $\mathbb{R}^2$. Let $f$ be a continuous function defined on $S=\mathbf{X}(D)$.

Suppose $D$ is partitioned into finitely many small pieces $D_1,D_2,\cdots,D_m$, every two of them can only intersect along some parts of their boundaries, and the area of the $k$-th pieces is $\triangle D_k$. $P=\{D_1,\cdots,D_m\} $. Let $\xi_k$ denote an arbitrary ``test point'' in $D_k$.

Denote $J = \{1,2, \cdots,m\}$. Then,

\begin{itemize}
\item[i)]  For a fixed natural number $K \in N^{+} (1 \leq K < m)$, denote $J_K = \{i_1 , \cdots ,i_K \} \subset J$, and $J\setminus J_{K}=\{1,2,\cdots,m\}\setminus \{i_1 , \cdots ,i_K \} $.
Then,

\begin{equation}
\begin{array}{ll}
\displaystyle \lim_{\max\{\triangle D_k\} \rightarrow 0} \sum\limits_{k\in J\setminus J_{K}} f(\mathbf{X}(\xi_k))\|\frac{\partial X}{\partial u} \times \frac{\partial X}{\partial v} \|_{\xi_k} \triangle D_k=\displaystyle \iint_{S} f dS.
\end{array}
\end{equation}

\item[ii)] If $P$ is an equal partition, for any natural number $K(m) (1 \leq K(m) < m)$, satisfying $\lim\limits_{m\rightarrow +\infty}\displaystyle \frac{K(m)}{m}=0 $, denote $J_{K(m)}=\{i_1 , \cdots ,i_{K(m)} \} \subset J$, and $J\setminus J_{K(m)}=\{1,2,\cdots,m\}\setminus \{i_1 , \cdots ,i_{K(m)} \} $.
Then,
\begin{equation}
\begin{array}{ll}
\displaystyle \lim_{\max\{\triangle D_k\} \rightarrow 0} \sum\limits_{k\in J\setminus J_{K(m)}} f(\mathbf{X}(\xi_k))\|\frac{\partial X}{\partial u} \times \frac{\partial X}{\partial v} \|_{\xi_k} \triangle D_k =\displaystyle \iint_{S} f dS.
\end{array}
\end{equation}

\end{itemize}
where $\mathbf{N}(u,v)=\displaystyle \frac{\partial \mathbf{X}}{\partial u} \times \frac{\partial \mathbf{X}}{\partial v}$.

\hspace*{1cm}

\textbf{Theorem 11} (Deleting Items Theorem of Scaler Surface Integral) Let $S$ be a smooth parameterized surface, that is $S$ is the image of the $C^1$ map $\mathbf{X}: D\rightarrow \mathbb{R}^3$,$(u,v)\in D$ , where $D$ is a connected, bounded region in $\mathbb{R}^2$. Let $f$ be a continuous function defined on $S=\mathbf{X}(D)$.

Suppose $D$ is partitioned into finitely many small pieces $D_1,D_2,\cdots,D_m$, every two of them can only intersect along some parts of their boundaries, and the area of the $k$-th pieces is $\triangle D_k$. $P=\{D_1,\cdots,D_m\} $. Let $\xi_k$ denote an arbitrary ``test point'' in $D_k$.

Let $\tilde{D}_1,\tilde{D}_2,\cdots,\tilde{D}_m$ be a distorted partition of $D_1,D_2,\cdots,D_m$, the area of the $D_k$ is $\triangle \tilde{D}_k$, $\tilde{P}=\{\tilde{D}_1,\cdots,\tilde{D}_m\}$, satisfying
\begin{equation*}
\displaystyle \lim\limits_{m\rightarrow +\infty} \sum\limits_{k=1}^{m} m(D_k \triangle \tilde{D}_k)=0
\end{equation*}
Then,
\begin{equation}
\begin{array}{ll}
\displaystyle \lim_{\max\{\triangle \tilde{D}_k\} \rightarrow 0} \sum\limits_{k=1}^m
f(\mathbf{X}(\xi_k))\|\frac{\partial X}{\partial u} \times \frac{\partial X}{\partial v} \|_{\xi_k} \triangle D_k=\displaystyle \iint_{S} f dS.
\end{array}
\end{equation}
where $\mathbf{N}(u,v)=\displaystyle \frac{\partial \mathbf{X}}{\partial u} \times \frac{\partial \mathbf{X}}{\partial v}$.

\hspace*{1cm}

\textbf{Theorem 12} Let $S$ be a smooth parameterized surface, that is $S$ is the image of the $C^1$ map $\mathbf{X}: D\rightarrow \mathbb{R}^3$, $(u,v)\in D$ , where $D$ is a connected, bounded region in $\mathbb{R}^2$. Let $f$ be a continuous function defined on $S=\mathbf{X}(D)$.

Suppose $D$ is partitioned into finitely many small pieces $D_1,D_2,\cdots,D_m$, every two of them can only intersect along some parts of their boundaries, and the area of the $k$-th pieces is $\triangle D_k$. $P=\{D_1,\cdots,D_m\} $. Let $\xi_k$ denote an arbitrary ``test point'' in $D_k$.

Let $\tilde{D}_1,\tilde{D}_2,\cdots,\tilde{D}_m$ be a distorted partition of $D_1,D_2,\cdots,D_m$, the area of the $D_k$ is $\triangle \tilde{D}_k$, $\tilde{P}=\{\tilde{D}_1,\cdots,\tilde{D}_m\}$, satisfying
\begin{equation*}
\displaystyle \lim\limits_{m\rightarrow +\infty} \sum\limits_{k=1}^{m} m(D_k \triangle \tilde{D}_k)=0
\end{equation*}
Denote $J = \{1,2, \cdots,m\}$. Then,
\begin{itemize}
\item[i)]  For a fixed natural number $K \in N^{+} (1 \leq K < m)$, denote $J_K = \{i_1 , \cdots ,i_K \} \subset J$, and $J\setminus J_{K}=\{1,2,\cdots,m\}\setminus \{i_1 , \cdots ,i_K \} $.
Then,

\begin{equation}
\begin{array}{ll}
\displaystyle \lim_{\max\{\triangle \tilde{D}_k\} \rightarrow 0} \sum\limits_{k\in J\setminus J_{K}} f(\mathbf{X}(\xi_k))\|\frac{\partial X}{\partial u} \times \frac{\partial X}{\partial v} \|_{\xi_k} \triangle D_k=\displaystyle \iint_{S} f dS.
\end{array}
\end{equation}

\item[ii)]If $P$ is an equal partition, for any natural number $K(m) (1 \leq K(m) < m)$, satisfying $\lim\limits_{m\rightarrow +\infty}\displaystyle \frac{K(m)}{m}=0 $, denote $J_{K(m)}=\{i_1 , \cdots ,i_{K(m)} \} \subset J$, and $J\setminus J_{K(m)}=\{1,2,\cdots,m\}\setminus \{i_1 , \cdots ,i_{K(m)} \} $.
Then,
\begin{equation}
\begin{array}{ll}
\displaystyle \lim_{\max\{\triangle \tilde{D}_k\} \rightarrow 0} \sum\limits_{k\in J\setminus J_{K(m)}} f(\mathbf{X}(\xi_k))\|\frac{\partial X}{\partial u} \times \frac{\partial X}{\partial v} \|_{\xi_k} \triangle D_k=\displaystyle \iint_{S} f dS.
\end{array}
\end{equation}
\end{itemize}
where $\mathbf{N}(u,v)=\displaystyle \frac{\partial \mathbf{X}}{\partial u} \times \frac{\partial \mathbf{X}}{\partial v}$.


The deleting items and disturbing mesh theorems of vector surface integral are as follows.

\hspace*{1cm}

\textbf{Theorem 13} (Deleting Items Theorem of Vector Surface Integral) Let $\mathbf{X}: D\rightarrow \mathbb{R}^3$ be a smooth parameterized surface, where $D$ is a bounded region in the $(u,v)$ plane, and let $\mathbf{F}(x,y,z)$ be a continuous vector field whose domain includes $S=\mathbf{X}(D)$.

Suppose $D$ is partitioned into finitely many small pieces $D_1,D_2,\cdots,D_m$, every two of them can only intersect along some parts of their boundaries, and the area of the $k$-th pieces is $\triangle D_k$. $P=\{D_1,\cdots,D_m\} $. Let $\xi_k$ denote an arbitrary ``test point'' in $D_k$.

Denote $J = \{1,2, \cdots,m\}$. Then,

\begin{itemize}
\item[i)]  For a fixed natural number $K \in N^{+} (1 \leq K < m)$, denote $J_K = \{i_1 , \cdots ,i_K \} \subset J$, and $J\setminus J_{K}=\{1,2,\cdots,m\}\setminus \{i_1 , \cdots ,i_K \} $.
Then,
\begin{equation}
\begin{array}{ll}
\displaystyle \lim_{\max\{\triangle D_k\} \rightarrow 0} \sum\limits_{k\in J\setminus J_{K}}
\mathbf{F}(\mathbf{X}(\xi_k))\cdot \mathbf{N}(\xi_k) \triangle D_k
=\displaystyle \iint_{\mathbf{X}} \mathbf{F}\cdot d\mathbf{S}.
\end{array}
\end{equation}

\item[ii)] If $P$ is an equal partition, for any natural number $K(m) (1 \leq K(m) < m)$, satisfying $\lim\limits_{m\rightarrow +\infty}\displaystyle \frac{K(m)}{m}=0 $, denote $J_{K(m)}=\{i_1 , \cdots ,i_{K(m)} \} \subset J$, and $J\setminus J_{K(m)}=\{1,2,\cdots,m\}\setminus \{i_1 , \cdots ,i_{K(m)} \} $.
Then,
\begin{equation}
\begin{array}{ll}
\displaystyle \lim_{\max\{\triangle D_k\} \rightarrow 0} \sum\limits_{k\in J\setminus J_{K(m)}} \mathbf{F}(\mathbf{X}(\xi_k))\cdot \mathbf{N}(\xi_k) \triangle D_k
=\displaystyle \iint_{\mathbf{X}} \mathbf{F}\cdot d\mathbf{S}.
\end{array}
\end{equation}
\end{itemize}
where $\mathbf{N}(u,v)=\displaystyle \frac{\partial \mathbf{X}}{\partial u} \times \frac{\partial \mathbf{X}}{\partial v}$.

\hspace*{1cm}

\textbf{Theorem 14} (Disturbing Mesh Theorem of Vector Surface Integral) Let $\mathbf{X}: D\rightarrow \mathbb{R}^3$ be a smooth parameterized surface, where $D$ is a bounded region in the $(u,v)$ plane, and let $\mathbf{F}(x,y,z)$ be a continuous vector field whose domain includes $S=\mathbf{X}(D)$.

Suppose $D$ is partitioned into finitely many small pieces $D_1,D_2,\cdots,D_m$, every two of them can only intersect along some parts of their boundaries, and the area of the $k$-th pieces is $\triangle D_k$. $P=\{D_1,\cdots,D_m\} $. Let $\xi_k$ denote an arbitrary ``test point'' in $D_k$.

Let $\tilde{D}_1,\tilde{D}_2,\cdots,\tilde{D}_m$ be a distorted partition of $D_1,D_2,\cdots,D_m$, the area of the $D_k$ is $\triangle \tilde{D}_k$, $\tilde{P}=\{\tilde{D}_1,\cdots,\tilde{D}_m\}$, satisfying
\begin{equation*}
\displaystyle \lim\limits_{m\rightarrow +\infty} \sum\limits_{k=1}^{m} m(D_k \triangle \tilde{D}_k)=0
\end{equation*}
Then,
\begin{equation}
\begin{array}{ll}
\displaystyle \lim_{\max\{\triangle \tilde{D}_k\} \rightarrow 0}
\sum\limits_{k=1}^m \mathbf{F}(\mathbf{X}(\xi_k))\cdot \mathbf{N}(\xi_k) D_k
=\displaystyle \iint_{\mathbf{X}} \mathbf{F}\cdot d\mathbf{S}.
\end{array}
\end{equation}
where $\mathbf{N}(u,v)=\displaystyle \frac{\partial \mathbf{X}}{\partial u} \times \frac{\partial \mathbf{X}}{\partial v}$.

\hspace*{1cm}

\textbf{Theorem 15} (Deleting Items and Disturbing Mesh Theorem of Vector Surface Integral) Let $\mathbf{X}: D\rightarrow \mathbb{R}^3$ be a smooth parameterized surface, where $D$ is a bounded region in the $(u,v)$ plane, and let $\mathbf{F}(x,y,z)$ be a continuous vector field whose domain includes $S=\mathbf{X}(D)$.

Suppose $D$ is partitioned into finitely many small pieces $D_1,D_2,\cdots,D_m$, every two of them can only intersect along some parts of their boundaries, and the area of the $k$-th pieces is $\triangle D_k$. $P=\{D_1,\cdots,D_m\} $. Let $\xi_k$ denote an arbitrary ``test point'' in $D_k$.

Let $\tilde{D}_1,\tilde{D}_2,\cdots,\tilde{D}_m$ be a distorted partition of $D_1,D_2,\cdots,D_m$, the area of the $D_k$ is $\triangle \tilde{D}_k$, $\tilde{P}=\{\tilde{D}_1,\cdots,\tilde{D}_m\}$, satisfying
\begin{equation*}
\displaystyle \lim\limits_{m\rightarrow +\infty} \sum\limits_{k=1}^{m} m(D_k \triangle \tilde{D}_k)=0
\end{equation*}
Denote $J = \{1,2, \cdots,m\}$. Then,
\begin{itemize}
\item[i)]  For a fixed natural number $K \in N^{+} (1 \leq K < m)$, denote $J_K = \{i_1 , \cdots ,i_K \} \subset J$, and
          $J\setminus J_{K}=\{1,2,\cdots,m\}\setminus \{i_1 , \cdots ,i_K \} $.
Then,
\begin{equation}
\begin{array}{ll}
\displaystyle \lim_{\max\{\triangle \tilde{D}_k\} \rightarrow 0} \sum\limits_{k\in J\setminus J_{K}}
\mathbf{F}(\mathbf{X}(\xi_k))\cdot \mathbf{N}(\xi_k) \triangle \tilde{D}_k
=\displaystyle \iint_{\mathbf{X}} \mathbf{F}\cdot d\mathbf{S}.
\end{array}
\end{equation}

\item[ii)] If $P$ is an equal partition, for any natural number $K(m) (1 \leq K(m) < m)$, satisfying $\lim\limits_{m\rightarrow +\infty}\displaystyle \frac{K(m)}{m}=0 $, denote $J_{K(m)}=\{i_1 , \cdots ,i_{K(m)} \} \subset J$, and $J\setminus J_{K(m)}=\{1,2,\cdots,m\}\setminus \{i_1 , \cdots ,i_{K(m)} \} $.
Then,
\begin{equation}
\begin{array}{ll}
\displaystyle \lim_{\max\{\triangle \tilde{D}_k\} \rightarrow 0} \sum\limits_{k\in J\setminus J_{K(m)}}
\mathbf{F}(\mathbf{X}(\xi_k))\cdot \mathbf{N}(\xi_k) \triangle \tilde{D}_k
=\displaystyle \iint_{\mathbf{X}} \mathbf{F}\cdot d\mathbf{S}.
\end{array}
\end{equation}
\end{itemize}
where $\mathbf{N}(u,v)=\displaystyle \frac{\partial \mathbf{X}}{\partial u} \times \frac{\partial \mathbf{X}}{\partial v}$.

Theorem 10-15 can be easily proved according to Theorem 1-3 of multiple integral, while the integral variables of the first type and second type surface integrals are parameterized variables $(u,v)$ in domain $D$. Thus, we develop the deleting item theorem, disturbing mesh theorem and deleting item and disturbing mesh theorem for multiple integral, line integral and surface integral respectively.

In addition, deleting item Riemann sum of multiple integral, line integral and surface integral can be treated as incomplete Riemann sum, but disturbing mesh Riemann sums of them are different from the traditional Riemann sum, we call them non-Riemann sums.

\section{Deleting Items and Disturbing Mesh Theorem of Green's Theorem, Gauss's Theorem and Stokes' Theorem}
\label{}

The famous three theorems Green's Theorem, Gauss's Theorem and Stokes' Theorem and relative definitions are mainly referred to [7].

\textbf{Definition 6} Let $f:D\subset \mathbb{R}^n\rightarrow \mathbb{R}$ be differentiable at $x=(x_1,x_2,\cdots,x_n)$. The gradient of $f$ at $x$ is
\begin{equation}
\nabla f =\displaystyle ( \frac{\partial f}{\partial x_1}, \frac{\partial f}{\partial x_2},\cdots,\frac{\partial f}{\partial x_n}).
\end{equation}
$\nabla f$ is a vector field and is called the gradient field of $f$.

\textbf{Definition 7} Let $\mathbf{F}:D\subset \mathbb{R}^n\rightarrow \mathbb{R}^n$ be a vector field of class $C^1$ on the open set $D$. the divergence of $\mathbf{F}=(f_1,f_2,\cdots,f_n)$ is the scalar function
\begin{equation}
Div \mathbf{F} =\displaystyle \sum_{j=1}^{n} \frac{\partial f_j}{\partial x_j} =\nabla \cdot \mathbf{F}.
\end{equation}

\textbf{Definition 8}. Let $\mathbf{F}:D\subset \mathbb{R}^3\rightarrow \mathbb{R}^3$, $\mathbf{F}=(f_1,f_2,f_3)$, be a vector field of class $C^1$ on the open set $D$. The curl (or rotor) of $\mathbf{F}$ is the vector field defined, formally, by the determinant
\begin{equation}
curl \mathbf{F} =
\left| \begin{array}{*{3}{ccc}}
 \vec{i}&\vec{j}&\vec{k}\\
\displaystyle \frac{\partial }{\partial x}&\displaystyle \frac{\partial }{\partial y}&\displaystyle \frac{\partial }{\partial z}\\
 f_1& f_2& f_3\\
\end{array}
\right|
=\nabla \times \mathbf{F}
\end{equation}
where $\nabla =\displaystyle ( \frac{\partial }{\partial x}, \frac{\partial }{\partial y},\frac{\partial }{\partial z})$.
\
\

\textbf{Green's theorem} Let $D$ be a closed, bounded region in $\mathbb{R}^2$ whose boundary $C=\partial D$ consists of finitely many simple, closed, piecewise $C^1$ curves. Orient the curves of $C$ so that $D$ is on the left as one traverses $C$. Let $\mathbf{F}(x,y)=\{P(x,y),Q(x,y)\}$ be a vector field of class $C^1$ throughout $D$. Then
\begin{equation}
\displaystyle \iint_{D} (\frac{\partial Q}{\partial x}-\frac{\partial P}{\partial y}) dxdy=\displaystyle \oint_{\partial D} Pdx+Qdy.
\end{equation}

\
\

\textbf{Gauss's Theorem} Let $V$ be a bounded solid region in $\mathbb{R}^3$ whose boundary $\partial V$ consists of finitely many simple, closed orientable surfaces, each of which is oriented by unit normals that point away from $V$ . Let $\mathbf{F}$ be a vector field of class $C^1$ whose domain includes $V$. Then
\begin{equation}
\displaystyle  \iiint_{V} \nabla \cdot \mathbf{F} dV=\displaystyle \varoiint_{\partial V} \mathbf{F} \cdot d \mathbf{S}.
\end{equation}

\
\

\textbf{Stokes' Theorem } Let $S$ be a bounded, piecewise smooth, oriented surface in $\mathbb{R}^{3}$. Suppose that $\partial S$ consists of finitely many piecewise $C^1$, simple, closed curves each of which is oriented consistently with $S$. Let $\mathbf{F}$ be a vector fields of class $C^1$ whose domain includes $S$. Then
\begin{equation}
\displaystyle \iint_{S} \nabla \times \mathbf{F} \cdot d\mathbf{S}=\displaystyle \oint_{\partial S} \mathbf{F} \cdot d\mathbf{s}.
\end{equation}

In the rest,  we will combine the deleting items theorem, disturbing mesh theorem and deleting items and disturbing mesh theorem in one theorem  for Green's Theorem,Gauss's Theorem and Stokes' Theorem  respectively.

\hspace*{1.5cm}

\textbf{Theorem 16} (Deleting Items and Disturbing Mesh Theorem of Green' Theorem) Let $D$ be a closed ,bounded region in $\mathbb{R}^2$ whose boundary $C=\partial D$ consists of finitely many simple, closed, piecewise $C^1$ curves. Orient the curves of $C$ so that $D$ is on the left as one traverses $C$. Let $\mathbf{F}(x,y)=\{P(x,y),Q(x,y)\}$ be a vector field of class $C^1$ throughout $D$.

Suppose $D$ is partitioned into finitely many small pieces $D_1,D_2,\cdots,D_m$, every two of them can only intersect along some parts of their boundaries, $D=\bigcup\limits_{k=1}^{m} D_k$, and the area of the $k$-th pieces is $\triangle D_k$, $P_D=\{D_1,\cdots,D_m\}$. Let $\xi_k$ denote an arbitrary ``test point'' in $D_k$.

Suppose $C=\partial D=\mathbf{x}(T)$, $t\in T\subset \mathbb{R}$, where $T$ is a union set of finite intervals determined by $C$. Let $t_0,t_1,t_2,\cdots,t_p$ be any partition points of $T$, denote $I_k=[t_{k-1},t_k]$, $T=\bigcup\limits_{k=1}^{p} I_k$. $\mathbf{x}(t_0)=\mathbf{x}(t_p)$ \footnote{Since $C$ is closed.}. $\triangle t_k=t_k-t_{k-1}$, $P_{\partial D}=\{I_1,\cdots,I_p\}$. Let $\eta_k$ denote an arbitrary ``test point'' in $I_k$.
\footnote{As long as $p$ is large enough, all intersect points of piecewise $C^1$ curves are adopted as points of partition $P$, hence, $C=\mathbf{x}(T)$ is $C^1$ curve in each $(t_{k-1},t_k)$.}

Let $\tilde{D}_1,\tilde{D}_2,\cdots,\tilde{D}_m$ be a distorted partition of $D_1,D_2,\cdots,D_m$, the area of $\tilde{D}_k$ is $\triangle \tilde{D}_k$, satisfying
\begin{equation*}
\displaystyle \lim\limits_{m\rightarrow +\infty} \sum\limits_{k=1}^{m} m(D_k \triangle \tilde{D}_k)=0
\end{equation*}

Let $\tilde{I}_1,\tilde{I}_2,\cdots,\tilde{I}_p$ be a distorted partition of $I_1,I_2,\cdots,I_p$,  the length of $\tilde{I}_k$ is $\triangle \tilde{t}_k$, satisfying
\begin{equation*}
\displaystyle \lim\limits_{p\rightarrow +\infty} \sum\limits_{k=1}^{p} m(I_k \triangle \tilde{I}_k)=0
\end{equation*}
Denote $J = \{1,2, \cdots,m\}$, $\hat{J} = \{1,2, \cdots,p\}$. Then,

\begin{itemize}
\item[i)]
\begin{equation}
\begin{array}{ll}
\displaystyle \lim_{\max\{\triangle \tilde{D}_k\} \rightarrow 0} \sum\limits_{k\in J}
\left| \begin{array}{rr}
\displaystyle \frac{\partial}{\partial x}  & \displaystyle \frac{\partial}{\partial y} \\
                           P & Q
\end{array}
\right|_{\xi_k}  \triangle \tilde{D}_k
=\displaystyle \iint_{D} (\frac{\partial Q}{\partial x}-\frac{\partial P}{\partial y}) dxdy \\
=\displaystyle \oint_{\partial D} Pdx+Qdy
=\displaystyle \lim_{\max\{\triangle \tilde{t}_k\} \rightarrow 0} \sum\limits_{k\in \hat{J}}
\mathbf{F}(\mathbf{x}(\eta_k)) \cdot \mathbf{x}^{'}(\eta_k) \triangle \tilde{t}_k
\end{array}
\end{equation}

\item[ii)]  For a fixed natural number $K \in N^{+} (1 \leq K < m)$, denote $J_K = \{i_1 , \cdots ,i_K \} \subset J$, and $J\setminus J_{K}=\{1,2,\cdots,m\}\setminus \{i_1 , \cdots ,i_K \} $.
    For a fixed natural number $L \in N^{+} (1 \leq L < p)$, denote $\hat{J}_L = \{j_1 , \cdots ,j_L \} \subset \hat{J}$, and $\hat{J} \setminus \hat{J}_{L}=\{1,2,\cdots,p\}\setminus \{j_1 , \cdots ,j_L \} $.
Then,

\begin{itemize}
\item[a)]
\begin{equation}
\begin{array}{ll}
\displaystyle \lim_{\max\{\triangle D_k\} \rightarrow 0} \sum\limits_{k\in J\setminus J_{K}}
\left| \begin{array}{rr}
\displaystyle \frac{\partial}{\partial x}  & \displaystyle \frac{\partial}{\partial y} \\
                           P & Q
\end{array}
\right|_{\xi_k}  \triangle D_k
=\displaystyle \iint_{D} (\frac{\partial Q}{\partial x}-\frac{\partial P}{\partial y}) dxdy \\
=\displaystyle \oint_{\partial D} Pdx+Qdy
=\displaystyle \lim_{\max\{\triangle t_k\} \rightarrow 0} \sum\limits_{k\in \hat{J}\setminus \hat{J}_{L}}
\mathbf{F}(\mathbf{x}(\eta_k)) \cdot \mathbf{x}^{'}(\eta_k) \triangle t_k
\end{array}
\end{equation}

\item[b)]
\begin{equation}
\begin{array}{ll}
\displaystyle \lim_{\max\{\triangle \tilde{D}_k\} \rightarrow 0} \sum\limits_{k\in J\setminus J_{K}}
\left| \begin{array}{rr}
\displaystyle \frac{\partial}{\partial x}  & \displaystyle \frac{\partial}{\partial y} \\
                           P & Q
\end{array}
\right|_{\xi_k}  \triangle \tilde{D}_k
=\displaystyle \iint_{D} (\frac{\partial Q}{\partial x}-\frac{\partial P}{\partial y}) dxdy \\
=\displaystyle \oint_{\partial D} Pdx+Qdy
=\displaystyle \lim_{\max\{\triangle \tilde{t}_k\} \rightarrow 0} \sum\limits_{k\in \hat{J}\setminus \hat{J}_{L}}
\mathbf{F}(\mathbf{x}(\eta_k)) \cdot \mathbf{x}^{'}(\eta_k) \triangle \tilde{t}_k
\end{array}
\end{equation}
\end{itemize}

\item[iii)]
 If both $P_D$ and $P_{\partial D}$ are equal partitions, for any natural number $K(m) (1 \leq K(m) < m)$, satisfying $\lim\limits_{m\rightarrow +\infty}\displaystyle \frac{K(m)}{m}=0 $, denote $J_{K(m)}=\{i_1 , \cdots ,i_{K(m)} \} \subset J$, and $J\setminus J_{K(m)}=\{1,2,\cdots,m\}\setminus \{i_1 , \cdots ,i_{K(m)} \} $. \\
    For any natural number $L(p) (1 \leq L(p) < p)$, satisfying $\lim\limits_{p\rightarrow +\infty}\displaystyle \frac{L(p)}{p}=0 $, denote $\hat{J}_{L(p)}=\{j_1 , \cdots ,j_{L(p)} \} \subset \hat{J}$, and $\hat{J}\setminus \hat{J}_{L(p)}=\{1,2,\cdots,p\}\setminus \{j_1 , \cdots ,j_{L(p)} \} $.
Then,

\begin{itemize}
\item[a)]
\begin{equation}
\begin{array}{ll}
\displaystyle \lim_{\max\{\triangle D_k\} \rightarrow 0} \sum\limits_{k\in J\setminus J_{K(m)}}
\left| \begin{array}{rr}
\displaystyle \frac{\partial}{\partial x}  & \displaystyle \frac{\partial}{\partial y} \\
                           P & Q
\end{array}
\right|_{\xi_k}  \triangle D_k
=\displaystyle \iint_{D} (\frac{\partial Q}{\partial x}-\frac{\partial P}{\partial y}) dxdy \\
=\displaystyle \oint_{\partial D} Pdx+Qdy
=\displaystyle \lim_{\max\{\triangle t_k\} \rightarrow 0} \sum\limits_{k\in \hat{J}\setminus \hat{J}_{L(p)}}
\mathbf{F}(\mathbf{x}(\eta_k)) \cdot \mathbf{x}^{'}(\eta_k) \triangle t_k
\end{array}
\end{equation}

\item[b)]
\begin{equation}
\begin{array}{ll}
\displaystyle \lim_{\max\{\triangle \tilde{D}_k\} \rightarrow 0} \sum\limits_{k\in J\setminus J_{K(m)}}
\left| \begin{array}{rr}
\displaystyle \frac{\partial}{\partial x}  & \displaystyle \frac{\partial}{\partial y} \\
                           P & Q
\end{array}
\right|_{\xi_k}  \triangle \tilde{D}_k
=\displaystyle \iint_{D} (\frac{\partial Q}{\partial x}-\frac{\partial P}{\partial y}) dxdy \\
=\displaystyle \oint_{\partial D} Pdx+Qdy
=\displaystyle \lim_{\max\{\triangle \tilde{t}_k\} \rightarrow 0} \sum\limits_{k\in \hat{J}\setminus \hat{J}_{L(p)}}
\mathbf{F}(\mathbf{x}(\eta_k)) \cdot \mathbf{x}^{'}(\eta_k) \triangle \tilde{t}_k
\end{array}
\end{equation}
\end{itemize}
where $d \mathbf{s} = \mathbf{x}^{'}(t)dt $, $\mathbf{x}^{'}(t)=\{x^{'}(t),y^{'}(t)\}$.
\end{itemize}

\hspace*{1cm}

\textbf{Theorem 17} (Deleting Items and Disturbing Mesh Theorem of Gauss' Theorem) Let $V$ be a bounded solid region in $\mathbb{R}^3$ whose boundary $\partial V$ consists of finitely many simple, closed orientable surfaces, each of which is oriented by unit normals that point away from $V$ . Let $\mathbf{F}$ be a vector field of class $C^1$ whose domain includes $V$.

Suppose $V$ is partitioned into finitely many small pieces $V_1,V_2,\cdots,V_m$, every two of them can only intersect along some parts of their boundaries, $D=\bigcup\limits_{k=1}^{m} V_k$, and the volume of the $V_k$ is $\triangle V_k$, $P_V= \{V_1,V_2,\cdots, V_m\}$. Let $\xi_k$ denote an arbitrary ``test point'' in $V_k$.

Suppose $S=\partial V=\mathbf{X}(D)$, $(u,v)\in D$, where $D$ are union of finite bounded sets determined by $S$. Suppose $D$ is partitioned into finitely many small pieces $D_1,D_2,\cdots,D_p$, every two of them can only intersect along some parts of their boundaries, $D=\bigcup\limits_{k=1}^{p} D_k$, and the area of the $D_k$ is $\triangle D_k$, $P_{\partial V}=\{D_1,D_2,\cdots,D_p\}$. Let $\eta_k$ denote an arbitrary ``test point'' in $D_k$.

Let $\tilde{V}_1,\tilde{V}_2,\cdots,\tilde{V}_m$ be a distorted partition of $V_1,V_2,\cdots,V_m$, the volume of $\tilde{V}_k$ is $\triangle \tilde{V}_k$, satisfying
\begin{equation*}
\displaystyle \lim\limits_{m\rightarrow +\infty} \sum\limits_{k=1}^{m} m(V_k\triangle \tilde{V}_k)=0
\end{equation*}

Let $\tilde{D}_1,\tilde{D}_2,\cdots,\tilde{D}_p$ be a distorted partition of $D_1,D_2,\cdots,D_p$, the area of $\tilde{D}_k$ is $\triangle \tilde{D}_k$,  satisfying
\begin{equation*}
\displaystyle \lim\limits_{p\rightarrow +\infty} \sum\limits_{k=1}^{p} m(D_k\triangle \tilde{D}_k)=0
\end{equation*}
Denote $J = \{1,2, \cdots,m\}$, $\hat{J} = \{1,2, \cdots,p\}$. Then,

\begin{itemize}
\item[i)]
\begin{equation}
\begin{array}{ll}
\displaystyle \lim_{\max\{\triangle \tilde{V}_k\} \rightarrow 0} \sum\limits_{k\in J}
\nabla \cdot \mathbf{F} (\xi_k) \triangle \tilde{V}_k
=\displaystyle  \iiint_{V}  \nabla \cdot \mathbf{F} d V \\
=\displaystyle \varoiint_{\partial V} \mathbf{F} \cdot d \mathbf{S}
=\displaystyle \lim_{\max\{\triangle \tilde{D}_k\} \rightarrow 0} \sum\limits_{k\in \hat{J}}
\mathbf{F}(\mathbf{X}(\eta_k))\cdot \mathbf{N}(\eta_k) \triangle \tilde{D}_k
\end{array}
\end{equation}

\item[ii)]
For a fixed natural number $K \in N^{+} (1 \leq K < m)$, denote $J_K = \{i_1 , \cdots ,i_K \} \subset J$, and $J\setminus J_{K}=\{1,2,\cdots,m\}\setminus \{i_1 , \cdots ,i_K \} $.
    For a fixed natural number $L \in N^{+} (1 \leq L < p)$, denote $\hat{J}_L = \{j_1 , \cdots ,j_L \} \subset \hat{J}$, and $\hat{J} \setminus \hat{J}_{L}=\{1,2,\cdots,p\}\setminus \{j_1 , \cdots ,j_L \} $.
Then,
\begin{itemize}
\item[a)]
\begin{equation}
\begin{array}{ll}
\displaystyle \lim_{\max\{\triangle V_k\} \rightarrow 0} \sum\limits_{k\in J\setminus J_{K}} \nabla \cdot \mathbf{F} (\xi_k) \triangle V_k
=\displaystyle  \iiint_{V}  \nabla \cdot \mathbf{F} d V \\
=\displaystyle \varoiint_{\partial V} F \cdot d \mathbf{S}
=\displaystyle \lim_{\max\{\triangle D_k\} \rightarrow 0} \sum\limits_{k\in \hat{J}\setminus \hat{J}_{L}}
 \mathbf{F}(\mathbf{X}(\eta_k))\cdot \mathbf{N}(\eta_k) \triangle D_k
\end{array}
\end{equation}

\item[b)]
\begin{equation}
\begin{array}{ll}
\displaystyle \lim_{\max\{\triangle \tilde{V}_k\} \rightarrow 0} \sum\limits_{k\in J\setminus J_{K}}
\nabla \cdot \mathbf{F} (\xi_k) \triangle \tilde{V}_k
=\displaystyle  \iiint_{V}  \nabla \cdot \mathbf{F} d V \\
=\displaystyle \varoiint_{\partial V} F \cdot d \mathbf{S}
=\displaystyle \lim_{\max\{\triangle \tilde{D}_k\} \rightarrow 0} \sum\limits_{k\in \hat{J}\setminus \hat{J}_{L}}
\mathbf{F}(\mathbf{X}(\eta_k))\cdot \mathbf{N}(\eta_k) \triangle \tilde{D}_k
\end{array}
\end{equation}
\end{itemize}

\item[iii)]  If both $P_V$ and $P_{\partial V}$ are equal partitions, for any natural number $K(m) (1 \leq K(m) < m)$, satisfying $\lim\limits_{m\rightarrow +\infty}\displaystyle \frac{K(m)}{m}=0 $, denote $J_{K(m)}=\{i_1 , \cdots ,i_{K(m)} \} \subset J$, and $J\setminus J_{K(m)}=\{1,2,\cdots,m\}\setminus \{i_1 , \cdots ,i_{K(m)} \} $. \\
    For any natural number $L(p) (1 \leq L(p) < p)$, satisfying $\lim\limits_{p\rightarrow +\infty}\displaystyle \frac{L(p)}{p}=0 $, denote $\hat{J}_{L(p)}=\{j_1 , \cdots ,j_{L(p)} \} \subset \hat{J}$, and $\hat{J}\setminus \hat{J}_{L(p)}=\{1,2,\cdots,p\}\setminus \{j_1 , \cdots ,j_{L(p)} \} $.
Then,
\begin{itemize}
\item[a)]
\begin{equation}
\begin{array}{ll}
\displaystyle \lim_{\max\{\triangle V_k\} \rightarrow 0} \sum\limits_{k\in J\setminus J_{K(m)}}
\nabla \cdot \mathbf{F} (\xi_k) \triangle V_k
=\displaystyle  \iiint_{V}  \nabla \cdot \mathbf{F} d V \\
=\displaystyle \varoiint_{\partial V} \mathbf{F} \cdot d \mathbf{S}
=\displaystyle \lim_{\max\{\triangle D_k\} \rightarrow 0} \sum\limits_{k\in \hat{J}\setminus \hat{J}_{L(p)}}
\mathbf{F}(\mathbf{X}(\eta_k))\cdot \mathbf{N}(\eta_k) \triangle D_k
\end{array}
\end{equation}

\item[b)]
\begin{equation}
\begin{array}{ll}
\displaystyle \lim_{\max\{\triangle \tilde{V}_k\} \rightarrow 0} \sum\limits_{k\in J\setminus J_{K(m)}}
\nabla \cdot \mathbf{F} (\xi_k) \triangle \tilde{V}_k
=\displaystyle  \iiint_{V}  \nabla \cdot \mathbf{F} d V \\
=\displaystyle \varoiint_{\partial V} \mathbf{F} \cdot d \mathbf{S}
=\displaystyle \lim_{\max\{\triangle \tilde{D}_k\} \rightarrow 0} \sum\limits_{k\in \hat{J}\setminus \hat{J}_{L(p)}}
\mathbf{F}(\mathbf{X}(\eta_k))\cdot \mathbf{N}(\eta_k) \triangle \tilde{D}_k
\end{array}
\end{equation}
\end{itemize}
where $\mathbf{N}(u,v)=\displaystyle \frac{\partial \mathbf{X}}{\partial u} \times \frac{\partial \mathbf{X}}{\partial v}$,
$d \mathbf{S}=\mathbf{N}(u,v) dudv$.
\end{itemize}

\hspace*{1cm}

\textbf{Theorem 18} (Deleting Items and Disturbing Mesh Theorem of Stokes' Theorem) Let $S$ be a bounded, piecewise smooth, oriented surface in $\mathbb{R}^{3}$. Suppose that $\partial S$ consists of finitely many piecewise $C^1$, simple, closed curves each of which is oriented consistently with $S$. Let $\mathbf{F}$ be a vector fields of class $C^1$ whose domain includes $S$.

Suppose $S=\mathbf{X}(D)$, $(u,v)\in D\subset \mathbb{R}^2$, where $D$ are union of finite bounded sets determined by $S$. Suppose $D$ is partitioned into finitely many small pieces $D_1,D_2,\cdots,D_p$, every two of them can only intersect along some parts of their boundaries, $D=\bigcup\limits_{k=1}^{p} D_k$, and the area of the $D_k$ is $\triangle D_k$, $P_{D}=\{D_1,D_2,\cdots,D_p\}$. Let $\xi_k$ denote an arbitrary ``test point'' in $D_k$.

Suppose $\partial S=\mathbf{x}(T)$, $t\in T\subset \mathbb{R}$, where $T$ is a union set of finite intervals determined by $\partial S$. Let $t_0,t_1,t_2,\cdots,t_p$ be any partition points of $T$, denote $I_k=[t_{k-1},t_k]$, $T=\bigcup\limits_{k=1}^{p} I_k$, $\mathbf{x}(t_0)=\mathbf{x}(t_p)$,
 $\triangle t_k= t_k-t_{k-1}$, $P_{\partial D}=\{I_1,I_2,\cdots, I_p\}$. Let $\eta_k$ denote an arbitrary ``test point'' in $I_k$

Let $\tilde{D}_1,\tilde{D}_2,\cdots,\tilde{D}_m$ be a distorted partition of $D_1,D_2,\cdots,D_m$, the area of $\tilde{D}_k$ is $\triangle \tilde{D}_k$, satisfying
\begin{equation*}
\displaystyle \lim\limits_{m\rightarrow +\infty} \sum\limits_{k=1}^{m} m(D_k \triangle \tilde{D}_k)=0
\end{equation*}

Let $\tilde{I}_1,\tilde{I}_2,\cdots,\tilde{I}_p$ be a distorted partition of $I_1,I_2,\cdots,I_p$, the length of $\tilde{I}_k$ is $\triangle \tilde{t}_k$, satisfying
\begin{equation*}
\displaystyle \lim\limits_{p\rightarrow +\infty} \sum\limits_{k=1}^{p} m(I_k \triangle \tilde{I}_k)=0
\end{equation*}
Denote $J=\{1,2, \cdots,m\}$, $\hat{J} = \{1,2, \cdots,p\}$. Then,

\begin{itemize}
\item[i)]
\begin{equation}
\begin{array}{ll}
\displaystyle \lim_{\max\{\triangle \tilde{D}_k\} \rightarrow 0} \sum\limits_{k\in J}
\nabla \times \mathbf{F}(\xi_k) \cdot \mathbf{N}(\xi_k) \triangle \tilde{D}_k
=\displaystyle \iint_{S} \nabla \times \mathbf{F} \cdot d\mathbf{S}\\
=\displaystyle \oint_{\partial S} \mathbf{F} \cdot d\mathbf{s}
=\displaystyle \lim_{\max\{\triangle \tilde{t}_k\} \rightarrow 0} \sum\limits_{k\in \hat{J}}
\mathbf{F}(\mathbf{x}(\eta_k)) \cdot \mathbf{x}^{'}(\eta_k) \triangle \tilde{t}_k
\end{array}
\end{equation}

\item[ii)]
For a fixed natural number $K \in N^{+} (1 \leq K < m)$, denote $J_K = \{i_1 , \cdots ,i_K \} \subset J$, and $J\setminus J_{K}=\{1,2,\cdots,m\}\setminus \{i_1 , \cdots ,i_K \} $.
    For a fixed natural number $L \in N^{+} (1 \leq L < p)$, denote $\hat{J}_L = \{j_1 , \cdots ,j_L \} \subset \hat{J}$, and $\hat{J} \setminus \hat{J}_{L}=\{1,2,\cdots,p\}\setminus \{j_1 , \cdots ,j_L \} $.
Then,
\begin{itemize}
\item[a)]
\begin{equation}
\begin{array}{ll}
\displaystyle \lim_{\max\{\triangle D_k\} \rightarrow 0} \sum\limits_{k\in J\setminus J_{K}}
\nabla \times \mathbf{F}(\xi_k) \cdot \mathbf{N}(\xi_k) \triangle \tilde{D}_k
=\displaystyle \iint_{S} \nabla \times \mathbf{F} \cdot d\mathbf{S}\\
=\displaystyle \oint_{\partial S} \mathbf{F} \cdot d\mathbf{s}
=\displaystyle \lim_{\max\{\triangle t_k\} \rightarrow 0} \sum\limits_{k\in \hat{J}\setminus \hat{J}_{L}}
\mathbf{F}(\mathbf{x}(\eta_k)) \cdot \mathbf{x}^{'}(\eta_k) \triangle t_k
\end{array}
\end{equation}

\item[b)]
\begin{equation}
\begin{array}{ll}
\displaystyle \lim_{\max\{\triangle \tilde{D}_k\} \rightarrow 0} \sum\limits_{k\in J\setminus J_{K}}
\nabla \times \mathbf{F}(\xi_k) \cdot \mathbf{N}(\xi_k) \triangle \tilde{D}_k
=\displaystyle \iint_{S} \nabla \times \mathbf{F} \cdot d\mathbf{S}\\
=\displaystyle \oint_{\partial S} \mathbf{F} \cdot d\mathbf{s}
=\displaystyle \lim_{\max\{\triangle \tilde{t}_k\} \rightarrow 0} \sum\limits_{k\in \hat{J}\setminus \hat{J}_{L}}
\mathbf{F}(\mathbf{x}(\eta_k)) \cdot \mathbf{x}^{'}(\eta_k) \triangle \tilde{t}_k
\end{array}
\end{equation}
\end{itemize}

\item[iii)]
 If both $P_D$ and $P_{\partial D}$ are equal partitions , for any natural number $K(m) (1 \leq K(m) < m)$, satisfying $\lim\limits_{m\rightarrow +\infty}\displaystyle \frac{K(m)}{m}=0 $, denote $J_{K(m)}=\{i_1 , \cdots ,i_{K(m)} \} \subset J$, and $J\setminus J_{K(m)}=\{1,2,\cdots,m\}\setminus \{i_1 , \cdots ,i_{K(m)} \} $. \\
    For any natural number $L(p) (1 \leq L(p) < p)$, satisfying $\lim\limits_{p\rightarrow +\infty}\displaystyle \frac{L(p)}{p}=0 $, denote $\hat{J}_{L(p)}=\{j_1 , \cdots ,j_{L(p)} \} \subset \hat{J}$, and $\hat{J}\setminus \hat{J}_{L(p)}=\{1,2,\cdots,p\}\setminus \{j_1 , \cdots ,j_{L(p)} \} $.
Then,
\begin{itemize}
\item[a)]
\begin{equation}
\begin{array}{ll}
\displaystyle \lim_{\max\{\triangle D_k\} \rightarrow 0} \sum\limits_{k\in J\setminus J_{K(m)}}
\nabla \times \mathbf{F}(\xi_k) \cdot \mathbf{N}(\xi_k) \triangle D_k
=\displaystyle \iint_{S} \nabla \times \mathbf{F} \cdot d\mathbf{S}\\
=\displaystyle \oint_{\partial S} \mathbf{F} \cdot d\mathbf{s}
=\displaystyle \lim_{\max\{\triangle t_k\} \rightarrow 0} \sum\limits_{k\in \hat{J}\setminus \hat{J}_{L(p)}}
\mathbf{F}(\mathbf{x}(\eta_k)) \cdot \mathbf{x}^{'}(\eta_k) \triangle t_k
\end{array}
\end{equation}

\item[b)]
\begin{equation}
\begin{array}{ll}
\displaystyle \lim_{\max\{\triangle \tilde{D}_k\} \rightarrow 0} \sum\limits_{k\in J\setminus J_{K(m)}}
\nabla \times \mathbf{F}(\xi_k) \cdot \mathbf{N}(\xi_k) \triangle \tilde{D}_k
=\displaystyle \iint_{S} \nabla \times \mathbf{F} \cdot d\mathbf{S}\\
=\displaystyle \oint_{\partial S} \mathbf{F} \cdot d\mathbf{s}
=\displaystyle \lim_{\max\{\triangle \tilde{t}_k\} \rightarrow 0} \sum\limits_{k\in \hat{J}\setminus \hat{J}_{L(p)}}
\mathbf{F}(\mathbf{x}(\eta_k)) \cdot \mathbf{x}^{'}(\eta_k) \triangle \tilde{t}_k
\end{array}
\end{equation}
\end{itemize}
where $\mathbf{N}(u,v)=\displaystyle \frac{\partial \mathbf{X}}{\partial u} \times \frac{\partial \mathbf{X}}{\partial v}$,
$d \mathbf{S}=\mathbf{N}(u,v) dudv$, $d \mathbf{s} = \mathbf{x}^{'}(t)dt $, $\mathbf{x}^{'}(t)=\{x^{'}(t),y^{'}(t)\}$.
\end{itemize}
\
\
According to Theorem 1-15  and conclusions of Green's Theorem, Gauss's Theorem and Stokes' Theorem , we can easily prove Theorem 16-18. It is well-known that the three classical theorems play important roles in natural science, such as  electrostatics, fluid dynamics, mathematics of physics and engineering,etc, Theorem 16-18 reveal that there are more different physical status and phenomena (incomplete Riemann Sum , non-Riemann Sum (especially, disturbing mesh Riemann Sum) ), which have different physical backgrounds and interpretation, also converging to the limit integral status of three classical theorems. Based on limit theory, numerous incomplete Riemann Sums and non-Riemann Sums incidental to a fixed Riemann Sum will lead to the complexity in cognition and mathematical modeling of the nature and objective world in theory and thought.

\section{Deleting Items and Disturbing Mesh Theorem of General Stokes' Theorem}
\label{}

In manifold analysis with $k-$form, exterior derivative, chart, $k-$dimensional surface, manifold and integral over form and orientable surface, Green's Theorem, Gauss's Theorem and Stokes' Theorem can be expressed in general Stokes' Theorem [2,3,6-11], which is also treated as an extension of high dimension form of line integral and surface integral. To avoid verbose repeating content of manifold analysis, we only list the main definitions involved in our discussion. Since the notations are slight different in literatures \footnote{For example, a chart $\{(A,\varphi)\}$ in [3] is different from that in [9]. In [3], $A$ denotes range or domain of action of the chart on the surface $S$. But in [9], $A$ denotes parameter domain.}, the notations involved in the discussion  mainly refer to [3].
\
\

There are two general Stokes' theorems in manifold analysis, one is on $k-$chain [11], or $k-$cube [9], or $k$-dimensional surface [11], the other is on manifold. Both are described as follows [3].

\hspace*{1cm}

\textbf{General Stokes' Theorem on $k-$dimensional surface}  If $S$ be an oriented piecewise smooth  $k-$dimensional compact surface with  boundary $\partial S$ in the domain $G\subset \mathbb{R}^{n}$, in which a smooth  $(k-1)$-form $\omega$ is defined. Then
\begin{equation}
\displaystyle \int_{S}  d\omega=\displaystyle \int_{\partial S}  \omega.
\end{equation}
in which the orientation of the boundary $\partial S$ is that induced by the orientation of $S$.

\hspace*{1cm}

\textbf{General Stokes' Theorem on Manifold} Let $M$ be an oriented smooth  $n$-dimensional manifold and $\omega$ a smooth differential form of degree  $n-1$ and compact support on  $M$. Then
\begin{equation}
\displaystyle \int_{M}  d\omega=\displaystyle \int_{\partial M}  \omega.
\end{equation}
where the orientation of the boundary $\partial M$ of the manifold $M$ is induced by the orientation of the manifold $M$. If $\partial M=\emptyset$, then $\displaystyle \int_{M}  d\omega =0$.

\hspace*{1cm}

According to the discussions of deleting items and disturbing mesh theorems on multiple integral, line integral, surface integral, Green's Theorem, Gauss' Theorem and Stokes' Theorem, we give the deleting items and disturbing mesh theorem of General Stokes' Theorem on $k-$dimensional surface and manifold as follows.  For convenience, we assume the union of any partition and its distorted partition still satisfies the theorem's precondition involved in discussion.

\hspace*{1cm}

\textbf{Theorem 19} (Deleting Items and Disturbing Mesh Theorem of General Stokes' Theorem on $k-$Dimensional Surface) If $S$ be an oriented piecewise smooth  $k$-dimensional compact surface with boundary $\partial S$ in the domain $G\subset \mathbb{R}^{n}$, in which a smooth  $(k-1)$-form $\omega$ is defined.

Suppose a finite atlas $\{(\tilde{U}_i,\varphi_i)\}_{i=1}^{r}$ covers $S$, where $S=\bigcup\limits_{i=1}^{r} \tilde{U}_i \subseteq \mathbb{R}^n$. There exists a partition $\{U_1,\cdots,U_m\}$ of $S$, $U_i\subseteq \tilde{U}_i\bigcap S$, $V=\bigcup\limits_{i=1}^{r} \varphi_i^{-1}(U_i)\subseteq \mathbb{R}^k$, and each two of $\{U_i\}_{i=1}^{r}$ or $\{\varphi_i^{-1}(U_i)\}_{i=1}^{r}$ can only intersect along some parts of their boundaries. Suppose $V$ is partitioned into finitely many small pieces $V_1,V_2,\cdots,V_m$ , every two of them can only intersect along some parts of their boundaries, $V=\bigcup\limits_{i=1}^{m} V_i$, and the area (volume) of $V_i$ is $\triangle V_i$, $P_S=\{V_1,V_2,\cdots,V_m\}$.

Suppose $\partial S$ is covered by a finite atlas $\{(\tilde{A}_j,\psi_j)\}_{j=1}^{q}$, where $\partial S=\bigcup\limits_{j=1}^{q} \tilde{A}_j \subseteq \mathbb{R}^n$.  There exists a partition $\{A_1,\cdots,A_q\}$ of $\partial S$, $A_i\subseteq \tilde{A}_i\bigcap \partial S$, $D=\bigcup\limits_{j=1}^{q} \psi_j^{-1}(A_j)\subseteq \mathbb{R}^{k-1}$, and each two of $\{A_j\}_{j=1}^{q}$ or $\{\psi_j^{-1}(A_j)\}_{j=1}^{q}$ can only intersect along some parts of their boundaries. Suppose $D$ is partitioned into finitely many small pieces $D_1,D_2,\cdots,D_p$ , every two of them can only intersect along some parts of their boundaries, $D=\bigcup\limits_{j=1}^{p} D_j$, and the area (volume) of $D_j$ is $\triangle D_j$, $P_{\partial S}=\{D_1,D_2,\cdots,D_p\}$

Let $\tilde{V}_1,\tilde{V}_2,\cdots,\tilde{V}_m$ be a distorted partition of $V_1,V_2,\cdots,V_m$, the area (volume) of $\tilde{V}_k$ is $\triangle \tilde{V}_k$, satisfying
\begin{equation*}
\displaystyle \lim\limits_{m\rightarrow +\infty} \sum\limits_{i=1}^{m} m( V_i \triangle \tilde{V}_i)=0
\end{equation*}

Let $\tilde{D}_1,\tilde{D}_2,\cdots,\tilde{D}_p$ be a distorted partition of $D_1,D_2,\cdots,D_p$, the area (volume) of $\tilde{D}_k$ is $\triangle \tilde{D}_k$, satisfying
\begin{equation*}
\displaystyle \lim\limits_{p\rightarrow +\infty} \sum\limits_{j=1}^{p} m(D_j\triangle \tilde{D}_j)=0
\end{equation*}
Denote $J = \{1,2, \cdots,m\}$, $\hat{J} = \{1,2, \cdots,p\}$. Then,

\begin{itemize}
\item[i)]
\begin{equation}
\begin{array}{ll}
\displaystyle \lim_{\max\{\triangle \tilde{V}_i\} \rightarrow 0} \sum\limits_{i\in J}
\sum\limits_{k=1}^{r} \int_{\varphi_k(\tilde{V}_i\bigcap \varphi_k^{-1}(U_k))} d\omega
=\displaystyle \int_{S}  d\omega \\
=\displaystyle \int_{\partial S}  \omega
=\displaystyle \lim_{\max\{\triangle \tilde{D}_j\} \rightarrow 0} \sum\limits_{j\in \hat{J}}
\sum\limits_{k=1}^{q}  \int_{\psi_k(\tilde{D}_j\bigcap \psi_k^{-1}(A_k))}  \omega
\end{array}
\end{equation}

\item[ii)]
For a fixed natural number $K \in N^{+} (1 \leq K < m)$, denote $J_K = \{i_1 , \cdots ,i_K \} \subset J$, and $J\setminus J_{K}=\{1,2,\cdots,m\}\setminus \{i_1 , \cdots ,i_K \} $.
    For a fixed natural number $L \in N^{+} (1 \leq L < p)$, denote $\hat{J}_L = \{j_1 , \cdots ,j_L \} \subset \hat{J}$, and $\hat{J} \setminus \hat{J}_{L}=\{1,2,\cdots,p\}\setminus \{j_1 , \cdots ,j_L \} $.
Then,
\begin{itemize}
\item[(a)]
\begin{equation}
\begin{array}{ll}
\displaystyle \lim_{\max\{\triangle V_i\} \rightarrow 0} \sum\limits_{i\in J\setminus J_{K}}
\sum\limits_{k=1}^{r} \int_{\varphi_k(V_i\bigcap \varphi_k^{-1}(U_k))} d\omega
=\displaystyle \int_{S}  d\omega \\
=\displaystyle \int_{\partial S}  \omega
=\displaystyle \lim_{\max\{\triangle D_j\} \rightarrow 0} \sum\limits_{j \in \hat{J}\setminus \hat{J}_{L}}
\sum\limits_{k=1}^{q}  \int_{\psi_k(D_j\bigcap \psi_k^{-1}(A_k))}  \omega
\end{array}
\end{equation}

\item[(b)]
\begin{equation}
\begin{array}{ll}
\displaystyle \lim_{\max\{\triangle \tilde{V}_i\} \rightarrow 0} \sum\limits_{i\in J\setminus J_{K}}
\sum\limits_{k=1}^{r} \int_{\varphi_k(\tilde{V}_i\bigcap \varphi_k^{-1}(U_k))} d\omega
=\displaystyle \int_{S}  d\omega  \\
=\displaystyle \int_{\partial S}  \omega
=\displaystyle \lim_{\max\{\triangle \tilde{D}_j\} \rightarrow 0} \sum\limits_{j \in \hat{J}\setminus \hat{J}_{L}}
\sum\limits_{k=1}^{q}  \int_{\psi_k(\tilde{D}_j\bigcap \psi_k^{-1}(A_k))}  \omega
\end{array}
\end{equation}
\end{itemize}

\item[iii)]
 If both $P_S$ and $P_{\partial S}$ are equal partitions, which make each $\varphi_i^{-1}(U_i)$ and $\psi_i^{-1}(A_j)$ into equal partition respectively, for any natural number $K(m) (1 \leq K(m) < m)$, satisfying $\lim\limits_{m\rightarrow +\infty}\displaystyle \frac{K(m)}{m}=0 $, denote $J_{K(m)}=\{i_1 , \cdots ,i_{K(m)} \} \subset J$, and $J\setminus J_{K(m)}=\{1,2,\cdots,m\}\setminus \{i_1 , \cdots ,i_{K(m)} \} $. \\
    For any natural number $L(p) (1 \leq L(p) < p)$, satisfying $\lim\limits_{p\rightarrow +\infty}\displaystyle \frac{L(p)}{p}=0 $, denote $\hat{J}_{L(p)}=\{j_1 , \cdots ,j_{L(p)} \} \subset \hat{J}$, and $\hat{J}\setminus \hat{J}_{L(p)}=\{1,2,\cdots,p\}\setminus \{j_1 , \cdots ,j_{L(p)} \} $.
Then,
\begin{itemize}
\item[(a)]
\begin{equation}
\begin{array}{ll}
\displaystyle \lim_{\max\{\triangle V_i\} \rightarrow 0} \sum\limits_{i\in J\setminus J_{K(m)}}
\sum\limits_{k=1}^{r} \int_{\varphi_k(V_i\bigcap \varphi_k^{-1}(U_k))} d\omega
=\displaystyle \int_{S}  d\omega  \\
=\displaystyle \int_{\partial S}  \omega
=\displaystyle \lim_{\max\{\triangle D_j\} \rightarrow 0} \sum\limits_{j\in \hat{J}\setminus \hat{J}_{L(p)}}
\sum\limits_{k=1}^{q}  \int_{\psi_k(D_j\bigcap \psi_k^{-1}(A_k))}  \omega
\end{array}
\end{equation}

\item[(b)]
\begin{equation}
\begin{array}{ll}
\displaystyle \lim_{\max\{\triangle \tilde{V}_i\} \rightarrow 0} \sum\limits_{i\in J\setminus J_{K(m)}}
\sum\limits_{k=1}^{r} \int_{\varphi_k(\tilde{V}_i \bigcap \varphi_k^{-1}(U_k))} d\omega
=\displaystyle \int_{S}  d\omega \\
=\displaystyle \int_{\partial S}  \omega
=\displaystyle \lim_{\max\{\triangle \tilde{D}_j\} \rightarrow 0} \sum\limits_{j\in \hat{J}\setminus \hat{J}_{L(p)}}
\sum\limits_{k=1}^{q}  \int_{\psi_k(\tilde{D}_j\bigcap \psi_k^{-1}(A_k))}  \omega
\end{array}
\end{equation}
\end{itemize}
\end{itemize}

\textbf{Proof} For partition $P_S=\{V_1,V_2,\cdots,V_m\}$,
\begin{equation*}
\begin{array}{ll}
\displaystyle S =\displaystyle  \bigcup \limits_{i=1}^{m} \{ \bigcup \limits_{k=1}^{r} \varphi_k({V_i \bigcap \varphi_k^{-1}(U_k)})\}
\end{array}
\end{equation*}

For partition $P_{\partial S}=\{D_1,D_2,\cdots,D_p\}$,
\begin{equation*}
\begin{array}{ll}
\displaystyle \partial S =\displaystyle  \bigcup \limits_{i=1}^{p} \{ \bigcup \limits_{k=1}^{q} \psi_k({D_i \bigcap \varphi_k^{-1}(A_k)})\}
\end{array}
\end{equation*}

(1) Since $ \displaystyle \lim\limits_{m\rightarrow +\infty} \sum\limits_{i=1}^{m} m(V_i\triangle \tilde{V}_i)=0 $,
according to the integral over $k$-form and $k$-dimensional surface, and Theorem 2, we obtain
\begin{equation*}
\begin{array}{ll}
\displaystyle \lim_{\max\{\triangle \tilde{V}_i\} \rightarrow 0} \sum\limits_{i\in J}
\int_{(\tilde{V}_i\triangle V_i)\bigcap \varphi_k^{-1}(U_k)} \varphi_k^{*}(d\omega) =0
\end{array}
\end{equation*}
\begin{equation*}
\begin{array}{ll}
\displaystyle \lim_{\max\{\triangle \tilde{V}_i\} \rightarrow 0} \sum\limits_{i\in J}
\sum\limits_{k=1}^{r} \int_{\varphi_k((\tilde{V}_i\triangle V_i)\bigcap \varphi_k^{-1}(U_k))} d\omega \\
=\displaystyle \lim_{\max\{\triangle \tilde{V}_i\} \rightarrow 0} \sum\limits_{i\in J}
\sum\limits_{k=1}^{r} \int_{(\tilde{V}_i\triangle V_i)\bigcap \varphi_k^{-1}(U_k)} \varphi_k^{*}(d\omega)=0
\end{array}
\end{equation*}
As $\bigcup\limits_{i=1}^{m} (\tilde{V}_i \setminus V_i) \subseteq  \bigcup\limits_{i=1}^{m} (\tilde{V}_i \triangle V_i)$,
we obtain,
\begin{equation*}
\begin{array}{ll}
\displaystyle \lim_{\max\{\triangle \tilde{V}_i\} \rightarrow 0} \sum\limits_{i\in J}
\sum\limits_{k=1}^{r} \int_{\varphi_k((\tilde{V}_i\setminus V_i)\bigcap \varphi_k^{-1}(U_k))} d\omega \\
=\displaystyle \lim_{\max\{\triangle \tilde{V}_i\} \rightarrow 0} \sum\limits_{i\in J}
\sum\limits_{k=1}^{r} \int_{(\tilde{V}_i\setminus V_i)\bigcap \varphi_k^{-1}(U_k)} \varphi_k^{*}(d\omega)=0
\end{array}
\end{equation*}
Then,
\begin{equation*}
\begin{array}{ll}
\displaystyle \lim_{\max\{\triangle \tilde{V}_i\} \rightarrow 0} \sum\limits_{i\in J}
\sum\limits_{k=1}^{r} \int_{\varphi_k(\tilde{V}_i\bigcap \varphi_k^{-1}(U_k))} d\omega
=\displaystyle \lim_{\max\{\triangle \tilde{V}_i\} \rightarrow 0} \sum\limits_{i\in J}
\sum\limits_{k=1}^{r} \int_{\tilde{V}_i\bigcap \varphi_k^{-1}(U_k)} \varphi_k^{*}(d\omega) \\
=\displaystyle \lim_{\max\{\triangle \tilde{V}_i\} \rightarrow 0} \sum\limits_{i\in J}
 \sum\limits_{k=1}^{r} \large(\int_{(\tilde{V}_i\setminus V_i)\bigcap \varphi_k^{-1}(U_k)} \varphi_k^{*}(d\omega)
 + \int_{ V_i\bigcap \varphi_k^{-1}(U_k)} \varphi_k^{*}(d\omega) \large) \\
=\displaystyle \lim_{\max\{\triangle \tilde{V}_i\} \rightarrow 0} \sum\limits_{i\in J}
 \sum\limits_{k=1}^{r} \int_{(\tilde{V}_i\setminus V_i)\bigcap \varphi_k^{-1}(U_k)} \varphi_k^{*}(d\omega)
 +\displaystyle \lim_{\max\{\triangle \tilde{V}_i\} \rightarrow 0} \sum\limits_{i\in J}
 \sum\limits_{k=1}^{r} \int_{ V_i\bigcap \varphi_k^{-1}(U_k)} \varphi_k^{*}(d\omega)  \\
=\displaystyle 0 + \lim_{\max\{\triangle V_i\} \rightarrow 0} \sum\limits_{i\in J}
 \sum\limits_{k=1}^{r} \int_{ V_i\bigcap \varphi_k^{-1}(U_k)} \varphi_k^{*}(d\omega)  \\
=\displaystyle \int_{S}  d\omega
\end{array}
\end{equation*}

Similar to above proof, we can obtain
\begin{equation*}
\displaystyle \lim_{\max\{\triangle \tilde{D}_j\} \rightarrow 0} \sum\limits_{j\in \hat{J}}
\sum\limits_{k=1}^{q}  \int_{\psi_k(\tilde{D}_j\bigcap \psi_k^{-1}(A_k))}  \omega
=\displaystyle \int_{\partial S}  \omega
\end{equation*}

According to general Stokes' Theorem on $k-$dimensional surface
\begin{equation*}
\begin{array}{ll}
\displaystyle \int_{S}  d\omega =\displaystyle \int_{\partial S}  \omega
\end{array}
\end{equation*}

Therefore, Theorem 19 i) holds.

(2) According to the definitions of integral over $k$-form and integral over $k$-dimensional surface and Theorem 1, $\int_{ V_i\bigcap \varphi_k^{-1}(U_k)} \varphi_k^{*}(d\omega)$ can be transformed into $k$-multiple integral. Then,
\begin{equation*}
\begin{array}{ll}
\displaystyle \lim_{\max\{\triangle V_i\} \rightarrow 0}
\sum\limits_{k=1}^{r} \int_{ V_i\bigcap \varphi_k^{-1}(U_k)} \varphi_k^{*}(d\omega)
=\displaystyle \sum\limits_{k=1}^{r} \lim_{\max\{\triangle V_i\} \rightarrow 0}
\int_{ V_i\bigcap \varphi_k^{-1}(U_k)} \varphi_k^{*}(d\omega)=0
\end{array}
\end{equation*}
Therefore,
\begin{equation}
\begin{array}{ll}
\displaystyle \int_{S}  d\omega=\displaystyle \lim_{\max\{\triangle V_i\} \rightarrow 0} \sum\limits_{i\in J}
\sum\limits_{k=1}^{r} \int_{ V_i\bigcap \varphi_k^{-1}(U_k)} \varphi_k^{*}(d\omega)\\
=\displaystyle \lim_{\max\{\triangle V_i\} \rightarrow 0} \sum\limits_{i\in J\setminus J_{K}}
\sum\limits_{k=1}^{r} \int_{ V_i\bigcap \varphi_k^{-1}(U_k)} \varphi_k^{*}(d\omega)
+\displaystyle \lim_{\max\{\triangle V_i\} \rightarrow 0} \sum\limits_{i\in J_{K}}
\sum\limits_{k=1}^{r} \int_{ V_i\bigcap \varphi_k^{-1}(U_k)} \varphi_k^{*}(d\omega)\\
=\displaystyle \lim_{\max\{\triangle V_i\} \rightarrow 0} \sum\limits_{i\in J\setminus J_{K}}
\sum\limits_{k=1}^{r} \int_{ V_i\bigcap \varphi_k^{-1}(U_k)} \varphi_k^{*}(d\omega)
+\displaystyle \sum\limits_{i\in J_{K}} \lim_{\max\{\triangle V_i\} \rightarrow 0}
\sum\limits_{k=1}^{r} \int_{ V_i\bigcap \varphi_k^{-1}(U_k)} \varphi_k^{*}(d\omega)\\
=\displaystyle \lim_{\max\{\triangle V_i\} \rightarrow 0} \sum\limits_{i\in J\setminus J_{K}}
\sum\limits_{k=1}^{r} \int_{ V_i\bigcap \varphi_k^{-1}(U_k)} \varphi_k^{*}(d\omega)
+0\\
=\displaystyle \lim_{\max\{\triangle V_i\} \rightarrow 0} \sum\limits_{i\in J\setminus J_{K}}
\sum\limits_{k=1}^{r} \int_{ V_i\bigcap \varphi_k^{-1}(U_k)} \varphi_k^{*}(d\omega)
\end{array}
\end{equation}

Again, as a chart belonging to an atlas of $S$ generates an atlas of $\partial S$ [3],  $\int_{D_j\bigcap \psi_k^{-1}(A_k)} \psi_k^{*}(\omega)$ is integral of $(k-1)$-form on $(k-1)$- dimensional surface, it can be transformed into $(k-1)$-multiple integral. Similar to the above proof, we can obtain
\begin{equation*}
\begin{array}{ll}
\displaystyle \lim_{\max\{\triangle D_j\} \rightarrow 0}
\sum\limits_{k=1}^{q}  \int_{\psi_k(D_j\bigcap \psi_k^{-1}(A_k))}  \omega
=\displaystyle \lim_{\max\{\triangle D_j\} \rightarrow 0}
\sum\limits_{k=1}^{q}  \int_{D_j\bigcap \psi_k^{-1}(A_k)}  \psi_k^{*}(\omega)=0
\end{array}
\end{equation*}
and,
\begin{equation*}
\begin{array}{ll}
\displaystyle \lim_{\max\{\triangle D_j\} \rightarrow 0} \sum\limits_{j \in \hat{J}\setminus \hat{J}_{L}}
\sum\limits_{k=1}^{q}  \int_{\psi_k(D_j\bigcap \psi_k^{-1}(A_k))}  \omega \\
=\displaystyle \lim_{\max\{\triangle D_j\} \rightarrow 0} \sum\limits_{j \in \hat{J}\setminus \hat{J}_{L}}
\sum\limits_{k=1}^{q}  \int_{D_j\bigcap \psi_k^{-1}(A_k)} \psi_k^{*}(\omega)
=\displaystyle \int_{\partial S} \omega
\end{array}
\end{equation*}

Again, according to general Stokes' theorem on $k$-dimensional surface and above proofs, Theorem 19 ii)a) holds.

(3) According to Theorem 19 i), Theorem 19 ii)a), and $\displaystyle \lim\limits_{m\rightarrow +\infty} \sum\limits_{i=1}^{m} m(V_i\triangle \tilde{V}_i)=0$, we obtain
\begin{equation*}
\begin{array}{ll}
\displaystyle \lim_{\max\{\triangle \tilde{V}_i\} \rightarrow 0} \sum\limits_{i\in J\setminus J_{K}}
\sum\limits_{k=1}^{r} \int_{\varphi_k(\tilde{V}_i\bigcap \varphi_k^{-1}(U_k))} d\omega\\
=\displaystyle \lim_{\max\{\triangle \tilde{V}_i\} \rightarrow 0} \sum\limits_{i\in J\setminus J_{K}}
\sum\limits_{k=1}^{r} \int_{\tilde{V}_i\bigcap \varphi_k^{-1}(U_k)} \varphi_k^{*}(d\omega) \\
=\displaystyle \sum\limits_{k=1}^{r} \lim_{\max\{\triangle \tilde{V}_i\} \rightarrow 0} \sum\limits_{i\in J\setminus J_{K}}
 \int_{\tilde{V}_i\bigcap \varphi_k^{-1}(U_k)} \varphi_k^{*}(d\omega) \\
=\displaystyle\sum\limits_{k=1}^{r} \lim_{\max\{\triangle V_i\} \rightarrow 0} \sum\limits_{i\in J\setminus J_{K}}
 \int_{V_i\bigcap \varphi_k^{-1}(U_k)} \varphi_k^{*}(d\omega) \\
=\displaystyle \lim_{\max\{\triangle V_i\} \rightarrow 0} \sum\limits_{i\in \setminus J_{K}}
\sum\limits_{k=1}^{r} \int_{\varphi_k(V_i\bigcap \varphi_k^{-1}(U_k))} d\omega\\
=\displaystyle \int_{S}  d\omega
\end{array}
\end{equation*}
Similarly, we obtain
\begin{equation*}
\begin{array}{ll}
\displaystyle \lim_{\max\{\triangle \tilde{D}_j\} \rightarrow 0} \sum\limits_{j \in \hat{J}\setminus \hat{J}_{L}}
\sum\limits_{k=1}^{q}  \int_{\psi_k(\tilde{D}_j\bigcap \psi_k^{-1}(A_k))}  \omega
=\displaystyle \int_{\partial S}  \omega
\end{array}
\end{equation*}
According to general Stokes' theorem on $k$-dimensional surface and above proofs, Theorem 19 ii)b) holds.

(4) Since $P_S$ and $P_{\partial S}$ are equal partitions, and each $\varphi_i^{-1}(U_i)$ and $\psi_i^{-1}(A_j)$ into equal partition respectively, according to Theorem 3, we obtain,
\begin{equation*}
\begin{array}{ll}
\displaystyle \lim_{\max\{\triangle V_i\} \rightarrow 0} \sum\limits_{i\in J\setminus J_{K(m)}}
\sum\limits_{k=1}^{r} \int_{\varphi_k(V_i\bigcap \varphi_k^{-1}(U_k))} d\omega\\
=\displaystyle \lim_{\max\{\triangle V_i\} \rightarrow 0} \sum\limits_{i\in J\setminus J_{K(m)}}
\sum\limits_{k=1}^{r} \int_{\varphi_k(V_i \bigcap \varphi_k^{-1}(U_k))} \varphi_k^{*}(d\omega)\\
=\displaystyle \sum\limits_{k=1}^{r} \lim_{\max\{\triangle V_i\} \rightarrow 0} \sum\limits_{i\in J\setminus J_{K(m)}}
 \int_{\varphi_k(V_i \bigcap \varphi_k^{-1}(U_k))} \varphi_k^{*}(d\omega)\\
=\displaystyle \sum\limits_{k=1}^{r} \lim_{\max\{\triangle V_i\} \rightarrow 0} \sum\limits_{i\in J}
 \int_{\varphi_k(V_i \bigcap \varphi_k^{-1}(U_k))} \varphi_k^{*}(d\omega)\\
=\displaystyle  \lim_{\max\{\triangle V_i\} \rightarrow 0} \sum\limits_{i\in J}
 \sum\limits_{k=1}^{r} \int_{\varphi_k(V_i \bigcap \varphi_k^{-1}(U_k))} \varphi_k^{*}(d\omega)\\
=\displaystyle \int_{S}  d\omega
\end{array}
\end{equation*}
Similarly, we obtain,
$\displaystyle \lim_{\max\{\triangle D_j\} \rightarrow 0} \sum\limits_{j\in \hat{J}\setminus \hat{J}_{L(p)}}
\sum\limits_{k=1}^{q}  \int_{\psi_k(D_j\bigcap \psi_k^{-1}(A_k))}  \omega=\displaystyle \int_{\partial S}\omega$

According to general Stokes' theorem on $k$-dimensional surface and above proofs, Theorem 19 iii)a) holds.

(5) According to Theorem 19 i) ,Theorem 19 iii)a) and Theorem 3, we obtain,

\begin{equation*}
\begin{array}{ll}
\displaystyle \lim_{\max\{\triangle \tilde{V}_i\} \rightarrow 0} \sum\limits_{i\in J\setminus J_{K(m)}}
\sum\limits_{k=1}^{r} \int_{\varphi_k(\tilde{V}_i\bigcap \varphi_k^{-1}(U_k))} d\omega\\
=\displaystyle \lim_{\max\{\triangle \tilde{V}_i\} \rightarrow 0} \sum\limits_{i\in J\setminus J_{K(m)}}
\sum\limits_{k=1}^{r} \int_{\varphi_k(\tilde{V}_i \bigcap \varphi_k^{-1}(U_k))} \varphi_k^{*}(d\omega)\\
=\displaystyle \sum\limits_{k=1}^{r} \lim_{\max\{\triangle \tilde{V}_i\} \rightarrow 0} \sum\limits_{i\in J\setminus J_{K(m)}}
 \int_{\varphi_k(\tilde{V}_i \bigcap \varphi_k^{-1}(U_k))} \varphi_k^{*}(d\omega)\\
=\displaystyle \sum\limits_{k=1}^{r} \lim_{\max\{\triangle V_i\} \rightarrow 0} \sum\limits_{i\in J\setminus J_{K(m)}}
 \int_{\varphi_k(V_i \bigcap \varphi_k^{-1}(U_k))} \varphi_k^{*}(d\omega)\\
=\displaystyle  \lim_{\max\{\triangle V_i\} \rightarrow 0} \sum\limits_{i\in J\setminus J_{K(m)}}
 \sum\limits_{k=1}^{r} \int_{\varphi_k(V_i \bigcap \varphi_k^{-1}(U_k))} \varphi_k^{*}(d\omega)\\
=\displaystyle \int_{S}  d\omega
\end{array}
\end{equation*}
Similarly, we obtain,
\begin{equation*}
\displaystyle \lim_{\max\{\triangle \tilde{D}_j\} \rightarrow 0} \sum\limits_{j\in \hat{J}\setminus \hat{J}_{L(p)}}\sum\limits_{k=1}^{q}  \int_{\psi_k(\tilde{D}_j\bigcap \psi_k^{-1}(A_k))}\omega =\displaystyle \int_{\partial S}\omega
\end{equation*}
According to general Stokes' theorem on $k$-dimensional surface and above proofs, Theorem 19 iii)b) holds.

Summarily,Theorem 19 holds.

\hfill $\Box$

\hspace*{1cm}

\textbf{Remark 1} Since $\{(\partial \tilde{U}_i,\varphi_i|_{\mathbb{R}^{k-1}})\}_{i=1}^{r}$ is an atlas of $\partial S$, we can adopt $\psi_i=\varphi_i|_{\mathbb{R}^{k-1}}$.

\textbf{Remark 2} If all $\varphi_i: I^k\rightarrow U\subset \mathbb{R}^n$, $\psi_j:I^{k-1}\rightarrow U\subset \mathbb{R}^n$, Theorem 19 will be more legible and can directly apply Theorem 1-3 to proof.

\hspace*{1cm}

\textbf{Theorem 20}(Deleting Items and Disturbing Mesh Theorem of General Stokes' Theorem on Manifold) Let $M$ be an oriented smooth  $n$-dimensional manifold and $\omega$ a smooth differential form of degree  $n-1$ and compact support on  $M$.

Suppose a finite atlas $\{(U_i,\varphi_i)\}_{i=1}^{r}$ covers $\mbox{supp}\ \omega \subseteq M$, where $\mbox{supp} \omega=\bigcup\limits_{i=1}^{r} U_i \subseteq \mathbb{R}^n$. $\varphi_i: B_i\rightarrow U_i$, $V=\bigcup\limits_{i=1}^{r} B_i$, and let $e_1,\cdots,e_r$ be a partition of unity subordinate to that covering $\mbox{supp} \omega$, and that $\mbox{supp} e_i\subset U_i$, $i=1,\cdots,r$. Suppose $V$ is partitioned into finitely many small pieces $V_1,V_2,\cdots,V_m$ , every two of them can only intersect along some parts of their boundaries, $V=\bigcup\limits_{i=1}^{m} V_i$, and the area (volume) of $V_i$ is $\triangle V_i$, $P_M=\{V_1,V_2,\cdots,V_m\}$. And, $\varphi_i^{*}(e_i \omega)$ is the coordinate representation of the form $e_i \omega |_{U_i}$ in the domain $B_i$ of variation of the coordinates of the corresponding local chart.

Suppose $\partial M$ is covered by a finite atlas $\{(\tilde{A}_j,\psi_j)\}_{j=1}^{q}$, where $\partial M=\bigcup\limits_{j=1}^{q} \tilde{A}_j \subseteq \mathbb{R}^n$.  There exists a partition $\{A_1,\cdots,A_q\}$ of $\partial M$, $A_i\subseteq \tilde{A}_i\bigcap \partial M$, $D=\bigcup\limits_{j=1}^{q} \psi_j^{-1}(A_j)\subseteq \mathbb{R}^{k-1}$, and each two of $\{A_j\}_{j=1}^{q}$ or $\{\psi_j^{-1}(A_j)\}_{j=1}^{q}$ can only intersect along some parts of their boundaries. Suppose $D$ is partitioned into finitely many small pieces $D_1,D_2,\cdots,D_p$ , every two of them can only intersect along some parts of their boundaries, $D=\bigcup\limits_{j=1}^{p} D_i$, and the area (volume) of $D_j$ is $\triangle D_j$, $P_{\partial M}=\{D_1,D_2,\cdots,D_p\}$

Let $\tilde{V}_1,\tilde{V}_2,\cdots,\tilde{V}_m$ be a distorted partition of $V_1,V_2,\cdots,V_m$, the area (volume) of $\tilde{V}_k$ is $\triangle \tilde{V}_k$, satisfying
\begin{equation*}
\displaystyle \lim\limits_{m\rightarrow +\infty} \sum\limits_{i=1}^{m} m( V_i \triangle \tilde{V}_i)=0
\end{equation*}

Let $\tilde{D}_1,\tilde{D}_2,\cdots,\tilde{D}_p$ be a distorted partition of $D_1,D_2,\cdots,D_p$, the area (volume) of $\tilde{D}_k$ is $\triangle \tilde{D}_k$, satisfying
\begin{equation*}
\displaystyle \lim\limits_{p\rightarrow +\infty} \sum\limits_{j=1}^{p} m(D_j\triangle \tilde{D}_j)=0
\end{equation*}
Denote $J = \{1,2, \cdots,m\}$, $\hat{J} = \{1,2, \cdots,p\}$. Then,

\begin{itemize}
\item[i)]
\begin{equation}
\begin{array}{ll}
\displaystyle \lim_{\max\{\triangle \tilde{V}_i\} \rightarrow 0} \sum\limits_{i\in J}
\sum\limits_{k=1}^{r} \int_{\tilde{V}_i \bigcap B_k} \varphi_k^{*}(e_k d\omega)
=\displaystyle \int_{M}  d\omega \\
=\displaystyle \int_{\partial M}  \omega
=\displaystyle \lim_{\max\{\triangle \tilde{D}_j\} \rightarrow 0} \sum\limits_{j\in \hat{J}}
\sum\limits_{k=1}^{q}  \int_{\psi_k(\tilde{D}_j\bigcap \psi_k^{-1}(A_k))}  \omega
\end{array}
\end{equation}

\item[ii)]
For a fixed natural number $K \in N^{+} (1 \leq K < m)$, denote $J_K = \{i_1 , \cdots ,i_K \} \subset J$, and $J\setminus J_{K}=\{1,2,\cdots,m\}\setminus \{i_1 , \cdots ,i_K \} $.
    For a fixed natural number $L \in N^{+} (1 \leq L < p)$, denote $\hat{J}_L = \{j_1 , \cdots ,j_L \} \subset \hat{J}$, and $\hat{J} \setminus \hat{J}_{L}=\{1,2,\cdots,p\}\setminus \{j_1 , \cdots ,j_L \} $.
Then,
\begin{itemize}
\item[(a)]
\begin{equation}
\begin{array}{ll}
\displaystyle \lim_{\max\{\triangle V_i\} \rightarrow 0} \sum\limits_{i\in J\setminus J_{K}}
\sum\limits_{k=1}^{r} \int_{V_i \bigcap B_k} \varphi_k^{*}(e_k d\omega)
=\displaystyle \int_{M}  d\omega \\
=\displaystyle \int_{\partial M}  \omega
=\displaystyle \lim_{\max\{\triangle D_j\} \rightarrow 0} \sum\limits_{j \in \hat{J}\setminus \hat{J}_{L}}
\sum\limits_{k=1}^{q}  \int_{\psi_k(D_j\bigcap \psi_k^{-1}(A_k))}  \omega
\end{array}
\end{equation}

\item[(b)]
\begin{equation}
\begin{array}{ll}
\displaystyle \lim_{\max\{\triangle \tilde{V}_i\} \rightarrow 0} \sum\limits_{i\in J\setminus J_{K}}
\sum\limits_{k=1}^{r} \int_{\tilde{V}_i \bigcap B_k} \varphi_k^{*}(e_k d\omega)
=\displaystyle \int_{M}  d\omega  \\
=\displaystyle \int_{\partial M}  \omega
=\displaystyle \lim_{\max\{\triangle \tilde{D}_j\} \rightarrow 0} \sum\limits_{j \in \hat{J}\setminus \hat{J}_{L}}
\sum\limits_{k=1}^{q}  \int_{\psi_k(\tilde{D}_j\bigcap \psi_k^{-1}(A_k))}  \omega
\end{array}
\end{equation}
\end{itemize}

\item[iii)]
 If both $P_M$ and $P_{\partial M}$ are equal partitions, which make each $\varphi_i^{-1}(U_i)$ and $\psi_i^{-1}(A_j)$ into equal partition respectively, for any natural number $K(m) (1 \leq K(m) < m)$, satisfying $\lim\limits_{m\rightarrow +\infty}\displaystyle \frac{K(m)}{m}=0 $, denote $J_{K(m)}=\{i_1 , \cdots ,i_{K(m)} \} \subset J$, and $J\setminus J_{K(m)}=\{1,2,\cdots,m\}\setminus \{i_1 , \cdots ,i_{K(m)} \} $. \\
    For any natural number $L(p) (1 \leq L(p) < p)$, satisfying $\lim\limits_{p\rightarrow +\infty}\displaystyle \frac{L(p)}{p}=0 $, denote $\hat{J}_{L(p)}=\{j_1 , \cdots ,j_{L(p)} \} \subset \hat{J}$, and $\hat{J}\setminus \hat{J}_{L(p)}=\{1,2,\cdots,p\}\setminus \{j_1 , \cdots ,j_{L(p)} \} $.
Then,
\begin{itemize}
\item[(a)]
\begin{equation}
\begin{array}{ll}
\displaystyle \lim_{\max\{\triangle V_i\} \rightarrow 0} \sum\limits_{i\in J\setminus J_{K(m)}}
\sum\limits_{k=1}^{r} \int_{V_i \bigcap B_k} \varphi_k^{*}(e_k d\omega)
=\displaystyle \int_{M}  d\omega  \\
=\displaystyle \int_{\partial M}  \omega
=\displaystyle \lim_{\max\{\triangle D_j\} \rightarrow 0} \sum\limits_{j\in \hat{J}\setminus \hat{J}_{L(p)}}
\sum\limits_{k=1}^{q}  \int_{\psi_k(D_j\bigcap \psi_k^{-1}(A_k))}  \omega
\end{array}
\end{equation}

\item[(b)]
\begin{equation}
\begin{array}{ll}
\displaystyle \lim_{\max\{\triangle \tilde{V}_i\} \rightarrow 0} \sum\limits_{i\in J\setminus J_{K(m)}}
\sum\limits_{k=1}^{r} \int_{\tilde{V}_i \bigcap B_k} \varphi_k^{*}(e_k d\omega)
=\displaystyle \int_{M}  d\omega \\
=\displaystyle \int_{\partial M}  \omega
=\displaystyle \lim_{\max\{\triangle \tilde{D}_j\} \rightarrow 0} \sum\limits_{j\in \hat{J}\setminus \hat{J}_{L(p)}}
\sum\limits_{k=1}^{q}  \int_{\psi_k(\tilde{D}_j\bigcap \psi_k^{-1}(A_k))}  \omega
\end{array}
\end{equation}
\end{itemize}
\end{itemize}

\textbf{Proof}
Since manifold $M$ is compact,  $\mbox{supp} \omega$ $\subset M$ is compact on manifold $M$, there exist partitions of Unity on $M$.
As $d \omega$ is $n$-form on $n$-dimensional manifold $M$, and $\omega$ is $n-1$-form on $n$-dimensional manifold $M$,
 $\displaystyle \int_{\Psi}  d\omega$ and $ \int_{\partial \Psi}  \omega$ will be finally changed to $k$-dimensional or $(k-1)$-dimensional multiple integral about parameter domain under the definition of integral of form.
Similar to the proof of Theorem 19, according to Theorem 1-3, integral on $n-1$-form and $n$-dimensional manifold , multiple integral of $n$-dimension and $n-1$-dimension, and the conclusion of general Stokes' Theorem, it is easy to prove that Theorem 20 holds. We omit the verbatim repetition of detail proof.

\hfill $\Box$

\textbf{Remark 3} Since $\{(\partial \tilde{U}_i,\varphi_i|_{\mathbb{R}^{n-1}})\}_{i=1}^{r}$ is an atlas of $\partial M$, we can adopt $\psi_i=\varphi_i|_{\mathbb{R}^{n-1}}$.

\textbf{Remark 4} Since a chart can be defined as $\varphi: \mathbb{R}^k\rightarrow U\subset \mathbb{R}^n$ or  $\varphi: [0,1]^k\rightarrow U\subset \mathbb{R}^n$, and both $S$ and $M$ are compact, we assume that all parameter domains for the charts involved in integral over chart and surface
are bounded in Theorem 19 and Theorem 20. The assumption is convenient for application of multiple integral with Riemann sum.

Therefore we extend deleting items and disturbing mesh theorems of general Stokes' Theorem on $k-$dimensinal surface and manifold.

\section{Conclusion}
\label{}

The deleting items  and disturbing mesh theorems are extended to multiple integral,line integral, surface integral, Green's Theorem, Gauss' Theorem ,  Stokes' Theorem and general Stokes' Theorem, which are the foundation of advanced calculus in modern mathematics. Our conclusions address the following mathematical philosophy thinking.

Multiple integral is foundation of measure theory of multivariate, hence deleting item and disturbing mesh theorem of multiple integral directly affects the measurement in high dimension space (including Lebesgue measure and  measure zero in the sense of Lebesgue, etc), calculation of probability distribution, moments of multivariate random variable and elements of stochastic process.

Incomplete Riemann sum and  non-Riemann sum constructed by deleting items and disturbing mesh in any partition have different background and meaning from traditional Riemann sum especially in Faraday's law and Amp$\grave{e}$re's law, however, their limit behaviors obey the same rule, which will provoke one's deep thought of mathematical nature for physical world. In some sense, the classical Green's theorem, Gauss' Theorem and Stokes' Theorem are rules hiding in more complicated mathematical sequences with different intuitive physical background. Existences of Incomplete Riemann sum and  non-Riemann sum also make us to rethink the meaning of grid design and approaching accuracy of integral and differential equation in computation simulation.

Deleting items and disturbing mesh theorems of general Stokes' theorem provide a deep view of integral on manifold and $k$-dimensional surface, so that we can scan the subtle configurations of limit,$k$-form and integral in high-dimension space.

\end{document}